\documentclass[final,11pt,3p]{elsarticle}

\usepackage{amssymb}
\usepackage{amsfonts}

\usepackage{amsmath}
\usepackage{amsthm}
\usepackage{enumerate}

\RequirePackage[OT1]{fontenc}
\RequirePackage{amsthm,amsmath, xcolor}
\RequirePackage[numbers]{natbib}
\RequirePackage[colorlinks,citecolor=blue,urlcolor=blue]{hyperref}
\usepackage{amscd,amsfonts,amssymb,latexsym,array,hhline,graphicx}
\usepackage{float}
\usepackage{lineno}

\newtheorem{theorem}{Theorem}[section]
\newtheorem{proposition}[theorem]{Proposition}

\newtheorem{lemma}[theorem]{Lemma}

\newtheorem{remark}{Remark}[section]
\newtheorem{corollary}[theorem]{Corollary}

\newcommand\cA{{\cal A}}
\newcommand\cC{{\cal C}}

\newcommand\cI{{\cal I}}
\newcommand\cG{{\cal G}}

\newcommand\cB{{\cal B}}

\newcommand\cN{{\cal N}}

\newcommand\cQ{{\cal Q}}

\newcommand\cR{{\cal R}}

\newcommand\ve{\varepsilon}

\newcommand\Er{\mbox{Err}}

\def\bbr{{\mathbb R}}

\def\text#1{\hbox{#1}}
\def\proof{{\noindent \bf Proof. }}
\def\endproof{\mbox{\ $\qed$}}
\def\ch{\mbox{ch}}

\def\E{{\bf E}}
\def\A{{\bf A}}
\def\P{{\bf P}}

\def\K{{\bf K}}
\def\V{{\bf V}}
\def\p{{\bf p}}
\def\g{{\bf g}}
\def\u{{\bf u}}
\def\l{{\bf l}}
\def\C{{\bf C}}
\def\D{{\bf D}}
\def\H{{\bf H}}

\def\G{{\bf G}}
\def\L{{\bf L}}
\def\U{{\bf U}}
\def\M{{\bf M}}
\def\m{{\bf m}}
\def\c{{\bf c}}
\def\l{{\bf l}}
\def\b{{\bf b}}

\def\r{{\bf r}}
\def\s{{\bf s}}
\def\t{{\bf t}}
\def\v{{\bf v}}
\def\q{{\bf q}}
\def\k{{\bf k}}

\def\Chi{{\bf 1}}

\def\d{\mathrm{d}}
\def\build #1_#2{\mathrel{\mathop{\kern 0pt #1}\limits_{#2}}}
\newcommand\tr{\mbox{tr}}
\newcommand\Trg{\mbox{Tr}}
\newcommand{\wh}{\widehat}
\newcommand{\wt}{\widetilde}
\newcommand{\zs}[1]{{\mathchoice{#1}{#1}{\lower.25ex\hbox{$\scriptstyle#1$}}
{\lower0.25ex\hbox{$\scriptscriptstyle#1$}}}}

\numberwithin{equation}{section}

\makeatletter
\def\ps@pprintTitle{%
 \let\@oddhead\@empty
 \let\@evenhead\@empty
 \def\@oddfoot{}%
 \let\@evenfoot\@oddfoot}
\makeatother

\journal{Journal of Multivariate Analysis}

\begin{document}

\begin{frontmatter}

\title{
Improved estimation via model selection method for semimartingale regressions based on discrete data
   }


%
%
%
%
%
%

\author[label1]{Evgeny Pchelintsev\corref{cor1}}
\cortext[cor1]{Corresponding author}
\ead{evgen-pch@yandex.ru}

\author[label1,label2]{Serguei Pergamenshchikov}
\ead{Serge.Pergamenchtchikov@univ-rouen.fr}

\author[label1]{Maria Povzun}
\ead{povzunyasha@gmail.com}

\address[label1]{International Laboratory of Statistics of Stochastic Processes and
Quantitative Finance, Tomsk State University, Tomsk, Russia}

\address[label2]{Laboratoire de Math\'ematiques Raphael Salem,
   Universit\'e de Rouen Normandie, Saint-Etienne-du-Rouvray, France}

\begin{abstract}
In this paper we study a high dimension (Big Data) semimartingale regression model observed in the discrete time moments. We study this model in a nonparametric  setting.
To this end  improved (shrinkage) estimation methods  are developed and the non-asymptotic comparison between shrinkage and least squares estimates are studied.  The nonparametric improvement effect for the shrinkage estimates showing the significant advantage with respect to the parametric case is established.
Then, an model selection method based on these estimates is developed.  Non-asymptotic sharp oracle inequalities for the constructed model selection procedure are obtained.
Constructive sufficient conditions for the observation frequency providing  the robust efficiency property in adaptive setting  are found.
As an example,  the regression model with non-Gaussian Ornstein--Uhlenbeck--L\'evy noises is considered. The results of Monte-Carlo simulations are given which confirm numerically the obtained theoretical results.
\end{abstract}





\end{frontmatter}

\section{Introduction}\label{sec:In}

\subsection{Problem and motivations}

In this paper we consider the linear regression model in continuous time
 \begin{equation}\label{sec:nparam-md-1}
  \d y_\zs{t}=\,\sum^{q}_\zs{j=1}\beta_\zs{j}\u_\zs{j}(t)\,\d t
 +\d \xi_\zs{t}\,,
  \quad 0\le t\le n\,,
 \end{equation}
where $(\u_\zs{j})_\zs{1\le j\le q}$ is a system of  known  linear independent
 $1$-periodic $\bbr\to\bbr$ functions, the noise process $(\xi_\zs{t})_\zs{t\ge 0}$ is an
unobservable semimartingale with unknown distribution.
The  process
\eqref{sec:nparam-md-1} is observed only at the discrete time moments
\begin{equation}
\label{sec:In.Obvs-1}
(y_\zs{t_j})_\zs{0\le j\le N}\,,
\quad t_\zs{j}=\frac{j}{p}
\quad\mbox{and}\quad
N=np
\,,
\end{equation}
where the observations frequency $p$ is some fixed integer number. We consider the model \eqref{sec:nparam-md-1} in the big data setting, i.e. under condition that
parameter dimension is more than number of observations, i.e. $q>N$.
Usually, in these cases
for statistical models with independent observations  one uses the following methods: Lasso algorithm (see, for example,
\cite{Tibshirani1996})
or the Dantzig selector method proposed in \cite{CandezTao2007}.
For dependent observations such models
 were studied, for example,   in \cite{DeGregorioIacus2012} and \cite{Fujimori2019}.
But in all these papers the number of parameters $q$ is known and, therefore, unfortunately,
these methods can't be used  to estimate, for example, the
number of parameters $q$.
 It should be noted, that the case when the number of parameters $q$ is unknown is one of challenging problems
  in the signal and image processing theory
(see, for example, \cite{BayisaZhouCronieYu2019, BeltaiefChernoyarovPergamenchtchikov2019} and the references therein).
 In this paper, similar to
 \cite{GaltchoukPergamenshchikov2019, GaltchoukPergamenshchikov2020},
 we study this problem  in  nonparametric setting which
 allows us to consider the  models \eqref{sec:nparam-md-1} with
 unknown $q$ or with $q=+\infty$.
So, we consider the following observations model
\begin{equation}\label{sec:In.1}
  \d y_t=S(t)\d t+\d \xi_\zs{t}\,,
  \quad 0\le t\le n\,,
 \end{equation}
where $S$ is an unknown $1$-periodic $\bbr\to\bbr$
function
from $\L_\zs{2}[0,1]$.

 Note that if $(\xi_\zs{t})_\zs{t\ge 0}$ is a Brownian motion, then we obtain the "signal plus white noise"  model (see, for example, \cite{IbragimovKhasminskii1981, Kutoyants1984, Kutoyants1994} and etc.)
which is widely used in statistical radiophysics. In the case, when the the process $(\xi_\zs{t})_\zs{t\ge 0}$
is a  semimartingale, the models \eqref{sec:In.1} are
very popular in important practice problems such that, for example,
 signals and images processing \cite{BarbuBeltaiefPergamenshchikov2018, BeltaiefChernoyarovPergamenchtchikov2019, Kassam1988},
  finance and insurance (see, for example,
\cite{BarndorfNielsenShephard2001, BerdjanePergamenchtchikov2013, KabanovPergamenshchikov2020} and the references therein).
For the first time the nonparametric regression models \eqref{sec:In.1} with general semimartingale noises
were introduced in \cite{KonevPergamenshchikov2009a, KonevPergamenshchikov2009b}. Here the
 adaptive efficient robust estimation methods based on the model selection approach are developed.
More precisely, to take into account the dependent observations in the framework of the models
\eqref{sec:In.1} in the papers
  \cite{HopfnerKutoyants2009, KonevPergamenshchikov2003, KonevPergamenshchikov2010}
it was used the Ornstein--Uhlenbeck noise processes, so called color Gaussian noises.
To consider  non-Gaussian observations
  it was introduced
  in the papers
  \cite{BarbuBeltaiefPergamenshchikov2018, KonevPergamenshchikov2012,
  KonevPergamenshchikov2015, KonevPergamenshchikovPchelintsev2014, Pchelintsev2013}
    impulse noises defined through the semi-Markov or L\'evy processes.
 Moreover, for similar models in the papers
\cite{PchelintsevPchelintsevPergamenshchikov2018a,
PchelintsevPchelintsevPergamenshchikov2018b, PchelintsevPergamenshchikov2018a, PchelintsevPergamenshchikov2018b,  PovzunPchelintsev2017}
 improved (shrinkage) nonparametric estimation procedures have been developed that significantly improve the
non asymptotic estimation quality  compared with usual methods.
  It should be emphasized, that in all these papers the improved estimation problems are studied only for the complete observations cases, i.e. when the all trajectory $(y_\zs{t})_\zs{0\le t\le n}$ is available for the observations. Therefore, these results can't be apply to study the model \eqref{sec:nparam-md-1}
in big data setting.

\subsection{Main contributions}
  Our main goal in this paper is to develop improved estimation methods
  for the model \eqref{sec:nparam-md-1} on the basis of the observations \eqref{sec:In.Obvs-1}.
  As an example for $(\xi_\zs{t})_\zs{t\ge 0}$, we consider the noise defined by non-Gaussian Ornstein--Uhlenbeck process with unknown distribution.
In the general framework of the observations model \eqref{sec:In.1} we develop special model selection methods based on the improved weighted least squares estimates.
  We recall, that for the first time such approach was proposed in \cite{FourdrinierPergamenshchikov2007} for spherical symmetric regression models in discrete time and in
   \cite{KonevPergamenshchikov2010} for Gaussian regression models in continuous time. It should be noted that for
the non spherically symmetric distributions    we can not use directly the well-known improved (shrinkage) estimators proposed in \cite{JamesStein1961}.  To apply the improved estimation methods to the general regression models in continuous time one needs to use the modifications of the James -- Stein shrinkage procedure proposed in \cite{KonevPergamenshchikovPchelintsev2014, Pchelintsev2013} for parametric models and developed in
     \cite{PchelintsevPchelintsevPergamenshchikov2018a,
PchelintsevPchelintsevPergamenshchikov2018b, PchelintsevPergamenshchikov2018a, PchelintsevPergamenshchikov2018b, PovzunPchelintsev2017}  for nonparametric ones.
Next, in this paper we study the estimation accuracy improvement effect based on discrete data \eqref{sec:In.Obvs-1}.
We show that the proposed improved estimators outperform the ordinary least squares estimates
 in non-asymptotic accuracy estimation uniformly over the observation frequency $p\ge 1$.
 It turns out, that in this case the estimation improvement  effect is the same as for the complete observations case, i.e. it is  established that uniformly over the observation frequency
the improvement of the non-asymptotic accuracy  is much more significant than for parametric models, since according to the well-known James - Stein formula, the accuracy improvement increases when the parametric dimension  increases (see, \cite{PchelintsevPchelintsevPergamenshchikov2018b}). Recall, that for the parametric models this dimension is always fixed, while for the nonparametric ones it tends to infinity, that is, it becomes arbitrarily large with an increase in the number of observations. Therefore, the gain from the application of improved methods is essentially increasing with respect to the parametric case. In the next step, we develop a special adaptive estimation tool
 through the model selection approach. More precisely, we find constructive sufficient conditions for the observation frequency and for the noise distributions under which
 through the robust risks
  we show sharp non-asymptotic oracle inequalities  for the proposed improved model selection procedures.
Finally, using these oracle inequalities, we provide the efficiency property for these procedures in adaptive setting.

\subsection{Plan of the paper}

The rest of the paper is organized as follows.
In Section ~\ref{sec:Mn-cnds} we give all necessary conditions which are used for the model \eqref{sec:In.1}.  In Section~\ref{sec:Imp} we construct the shrinkage
weighted least squares estimates and study the improvement effect. In Section~\ref{sec:Mo} we
develop improved model selection method. In Section \ref{sec:MainRs} we announce the main results of this paper.
A numerical Monte-Carlo  analysis is given in Section \ref{sec:Sim}. The main properties of the non-Gaussian Ornstein--Uhlenbeck processes are studied in
  Section \ref{sec:Main-prps}.  Section \ref{sec:proofs-1} contains the proofs of the all main results.
In the appendix (Section \ref{sec:Appendix}) it is postponed all auxiliary technical results.

\section{Main conditions}\label{sec:Mn-cnds}

First we assume that, the noise process $(\xi_\zs{t})_\zs{t\ge 0}$ in the model  \eqref{sec:In.1}
 is a square integrable semimartingale  with the values in the
Skorokhod space $\D[0,n]$
 such that, for any
function $f$ from $\L_\zs{2}[0,n]$, the stochastic integral
\begin{equation}\label{sec:In.2}
I_\zs{n}(f)=\int^n_\zs{0}f(s)\d \xi_\zs{s}
 \end{equation}
  has the following properties
\begin{equation}\label{sec:In.3}
\E_\zs{Q} I_\zs{n}(f)=0 \quad\mbox{and}\quad \E_\zs{Q}
I^2_\zs{n}(f)\le \varkappa_\zs{Q} \int^{n}_\zs{0}\,f^{2}(s)\d s
\,.
\end{equation}
 Here $\E_\zs{Q}$ denotes the expectation with respect to the distribution
 $Q$ of the noise process $(\xi_\zs{t})_\zs{0\le t\le n}$ on the space $\D[0,n]$,
 $\varkappa_\zs{Q}>0$ is
some positive constant depending on the distribution $Q$. As to the noise distribution $Q$ we assume that it is unknown and belongs to some
  family $\cQ_\zs{n}$ of probability distributions in $\D[0,n]$.
 For this problem we use the quadratic risk, which for any estimate $\wh{S}$, is defined as
\begin{equation}
\label{sec:In.4}
\cR_\zs{Q}(\wh{S},S):=
\E_\zs{Q,S}\,\|\wh{S}-S\|^2
\quad\mbox{and}\quad
\|f\|^2:=\int^1_\zs{0}f^2(t) \d t\,,
\end{equation}
where $\E_\zs{Q,S}$ stands for the expectation with respect to the
distribution $\P_\zs{Q,S}$ of the process \eqref{sec:In.1} with a fixed
distribution $Q$ of the noise $(\xi_\zs{t})_\zs{0\le t\le n}$
and a given function $S$. Moreover, in the case when the distribution
$Q$ is unknown
we use also  the robust risk
\begin{equation}\label{sec:In.6}
\cR^{*}(\wh{S},S)=\sup_\zs{Q\in\cQ_\zs{n}}\,
\cR_\zs{Q}(\wh{S},S)
\,.
\end{equation}
\noindent
For the distribution family  $\cQ_\zs{n}$ we assume the following condition

$\H_\zs{1})$ {\sl The family $\cQ_\zs{n}$
is such that
$$
\varkappa^{*}=\varkappa^{*}_\zs{n}=\sup_\zs{Q\in\cQ_\zs{n}}\,\varkappa_\zs{Q}
<\infty
$$
and for any $\epsilon>0$
$$
\lim_\zs{n\to\infty}\,n^{-\epsilon}\,\varkappa^{*}=0\,.
$$
}

\noindent
Now we consider the noise $(\xi_\zs{t})_\zs{t\ge 0}$ in
\eqref{sec:In.1} defined by a non-Gaussian Ornstein--Uhlenbeck process. Such processes are used in
the financial Black--Scholes type markets with  jumps (see, for
example,  \cite{BarndorfNielsenShephard2001}, and the references therein). Let the noise
process in \eqref{sec:In.1} obeys the equation
\begin{equation}\label{sec:Ex.1}
\d\xi_\zs{t} = a\xi_\zs{t}\d t+\d u_\zs{t}\,,\quad \xi_\zs{0}=0\,,
\end{equation}
where
$$
u_\zs{t} =
\varrho_\zs{1}\, w_\zs{t}+\varrho_\zs{2}\,z_\zs{t}
\quad\mbox{and}\quad
z_\zs{t}=x*(\mu-\wt{\mu})_\zs{t}\,.
$$
Here $(w_\zs{t})_\zs{t\ge 0}$ is
a standard Brownian motion,
 "$*$"\ denotes the stochastic integral with respect to the compensated jump measure $\mu(\d s\,,\,\d x)$ with  deterministic
compensator $\wt{\mu}(\d s\,\d x)=\d s\Pi(\d x)$, i.e.
$$
z_\zs{t}=\int_0^t\int_{\bbr_\zs{*}}v\,(\mu-\wt{\mu})(\d s \,\d v)
\quad\mbox{and}\quad
\bbr_\zs{*}=
\bbr\setminus \{0\}\,,
$$
\noindent
 $\Pi(\cdot)$ is the L\'evy measure on $\bbr_\zs{*}$, (see, for example in \cite{ContTankov2004}), such that
\begin{equation}\label{sec:Ex.1-00_mPi}
\Pi(x^{2})=1
\quad\mbox{and}\quad
\Pi(x^{8})
\,<\,\infty\,.
\end{equation}
We use the notation $\Pi(\vert x\vert^{m})=\int_\zs{\bbr_\zs{*}}\,\vert z\vert^{m}\,\Pi(\d z)$.
Moreover,  we assume that the nuisance parameters  $a\le 0$, $\varrho_\zs{1}$
and $\varrho_\zs{2}$ satisfy the conditions
\begin{equation}\label{sec:Ex.01-1}
-a_\zs{max}\le  a\le 0\,,\quad
0< \underline{\varrho}\le \varrho^{2}_\zs{1}
\quad\mbox{and}\quad
\sigma_\zs{Q}=\varrho^{2}_\zs{1}+\varrho^{2}_\zs{2}\,
\le
\varsigma^{*}
\,,
\end{equation}
where
 the bounds
$a_\zs{max}$, $\underline{\varrho}$ and $\varsigma^{*}$ are functions of $n$, i.e.
$a_\zs{max}=a_\zs{max}(n)$,
$\underline{\varrho}=\varrho_\zs{n}$
and $\varsigma^{*}=\varsigma^{*}_\zs{n}$, such that for any $\epsilon>0$
\begin{equation}\label{sec:Ex.01-2}
\lim_\zs{n\to\infty}\,\frac{a_\zs{max}(n)+\varsigma^{*}_\zs{n}}{n^{\epsilon}}=0
\quad\mbox{and}\quad
\liminf_\zs{n\to\infty}\,n^{\epsilon}\,
\underline{\varrho}_\zs{n}>0\,.
\end{equation}

\noindent
In this case $\cQ_\zs{n}$ is the family
of all distributions of process \eqref{sec:In.1} -- \eqref{sec:Ex.1} on the Skorokhod space
$\D[0,n]$ satisfying the conditions \eqref{sec:Ex.01-1} -- \eqref{sec:Ex.01-2}.
It should be noted that in view of Corollary~\ref{Co.sec:Stc.1} and the last inequality in \eqref{sec:Ex.01-1},
 the condition \eqref{sec:In.3} for the process \eqref{sec:Ex.1} holds with
$\varkappa_\zs{Q}=2\varsigma^{*}$.
Note also that the process
\eqref{sec:Ex.1} is conditionally-Gaussian square integrable semimartingale with respect to
$\sigma$-algebra $\cG=\sigma\{z_\zs{t}\,,\,t\ge 0\}$ which is generated by jump process $(z_t)_{t\ge 0}$.

\section{Improved estimation method}\label{sec:Imp}

For estimating the unknown function $S$ in \eqref{sec:In.1} we will use it's Fourier expansion
on the time grid $\{t_\zs{1},\ldots,t_\zs{p}\}$ defined in \eqref{sec:In.Obvs-1}. To this end we chose basic $1$ - periodic
functions $(\phi_\zs{j})_\zs{1\le j\le p}$  from $\L_\zs{2}[0,1]$, which are orthonormal on this grid, i.e.

\begin{equation}\label{innerprod}
(\phi_\zs{i},\phi_\zs{j})_p=\frac{1}{p}\sum_\zs{l=1}^{p}\, \phi_\zs{i}(t_\zs{l})\phi_\zs{j}(t_\zs{l})=\Chi_{\{i=j\}}\,.
\end{equation}
Note that we do not require the functions $(\phi_\zs{1\le j\le p})$ to be orthogonal in $\L_\zs{2}[0,1]$, it suffices to have only the property \eqref{innerprod}.
Assume that these functions are  bounded, i.e.
for some  constant $\phi_\zs{*}\ge 1$, which may be depend on $n$,
\begin{equation}\label{sec:In.3-00}
\sup_\zs{0\le j\le n}\,\sup_\zs{0\le t\le 1}\vert\phi_\zs{j}(t)\vert\,
\le\,
\phi_\zs{*}
<\infty\,.
\end{equation}
\noindent
We can take, for example, the Demmler - Reinsh spline basis defined in \cite{DemmlerReinsh1975}.
In the sequel we will use the trigonometric functions $(\Trg_\zs{j})_\zs{j\ge 1}$.
We recall that  $\Trg_\zs{1}\equiv 1$ and for $j\ge 2$
\begin{equation}\label{sec:In.5}
 \Trg_\zs{j}(x)= \sqrt 2
\left\{
\begin{array}{c}
\cos(2\pi[j/2] x)\, \quad\mbox{for even}\quad j \,;\\[4mm]
\sin(2\pi[j/2] x)\quad\mbox{for odd}\quad j\,.
\end{array}
\right.
\end{equation}
Here $[a]$ denotes the integer part of $a$.  It should be noted
that the trigonometric basis possesses  the property \eqref{innerprod}
only in the case when $p$ is odd. It is clear that if $p$ is even we can always to reduce the observation to $p-1$.
Note now that for any $t\in \{t_\zs{1},\ldots,t_\zs{p}\}$
we can represent the unknown function  as
\begin{equation}\label{sec:Imp.2}
S(t)=\sum_{j=1}^p \theta_\zs{j}\phi_\zs{j}(t)
\quad\mbox{and}\quad
\theta_\zs{j}=(S,\phi_\zs{j})_\zs{p}
\,.
\end{equation}
Therefore, to estimate the function $S$ it suffices to estimate the coefficients $\theta_\zs{j}$.
To this end we note that
$$
\theta_\zs{j}=(S,\phi_\zs{j})_\zs{p}=
\sum^{p}_\zs{l=1}\,\phi_\zs{j}(t_\zs{l})S(t_\zs{l})(t_\zs{l}-t_\zs{l-1})\,,
$$
and it is nature to replace the difference
$S(t_\zs{l})(t_\zs{l}-t_\zs{l-1})$ by the observation from the model  \eqref{sec:In.1}
as $y_\zs{t_\zs{l}}-y_\zs{t_\zs{l-1}}$. It should be note here also that to obtain a "good"
estimator one needs to use all observations $(y_\zs{t_\zs{l}})_\zs{1\le l\le np}$. Since
the functions $\phi_\zs{j}$ and $S$ are $1$ periodic, i.e. for any integers $1\le l\le p$ and $1\le k\le n$
$$
\phi_\zs{j}(t_\zs{l+kp})=\phi_\zs{j}(t_\zs{l})
\quad\mbox{and}\quad
S(t_\zs{l+kp})=S(t_\zs{l})
\,,
$$
therefore, we can represent the coefficients in \eqref{sec:Imp.2} as
$$
\theta_\zs{j}=(S,\phi_\zs{j})_\zs{p}=
\frac{1}{n}\sum^{np}_\zs{l=1}\,\phi_\zs{j}(t_\zs{l})S(t_\zs{l})(t_\zs{l}-t_\zs{l-1})\,.
$$
Using here the observations
\eqref{sec:In.Obvs-1}, we define the estimator as

\begin{equation}\label{sec:Imp.3}
\wh{\theta}_\zs{j}=
\frac{1}{n}\sum^{np}_\zs{l=1}\,\phi_\zs{j}(t_\zs{l})\,(y_\zs{t_\zs{l}}- y_\zs{t_\zs{l-1}})
= \frac{1}{n}\int^n_\zs{0}\,\psi_{j}(t)\,\d
y_\zs{t}\,
\end{equation}
with
$$
\psi_{j}(t)=\sum_{k=1}^{np}\phi_\zs{j}(t_k)\Chi_\zs{(t_{k-1},t_k]}(t)\,.
$$
For any functions $f,g$ from $\L_\zs{2}[0,1]$ we denote by $(f,g)$ the inner product in $\L_\zs{2}[0,1]$.
The system of the functions $(\psi_\zs{j})_\zs{1\le j\le p}$
is orthonormal in $\L_\zs{2}[0,1]$, i.e.
$$
(\psi_\zs{i},\psi_\zs{j})=\int_0^1\psi_\zs{i}(t)\psi_\zs{j}(t) \d t
=(\phi_\zs{i},\phi_\zs{j})_p=\Chi_{\{i=j\}},
$$
and  the corresponding Fourier coefficients in $\L_\zs{2}[0,1]$  can be represented
as
\begin{equation}\label{discoeff}
\bar{\theta}_\zs{j}=(S,\psi_\zs{j})=\theta_\zs{j}+h_\zs{j}\,,
\end{equation}
where
\begin{equation*}
h_\zs{j}=h_\zs{j}(S)=\sum_{k=1}^{p}\int_{t_{k-1}}^{t_k}\phi_\zs{j}(t_k)\left( S(t)-S(t_k)\right) \,\d t.
\end{equation*}
Therefore, we can represent the Fourier coefficients estimators as
\begin{equation}\label{sec:Imp.4}
\wh{\theta}_\zs{j}=\bar{\theta}_\zs{j}+\frac{1}{\sqrt{n}}\xi_\zs{j}
\quad\mbox{and}\quad
 \xi_\zs{j}=\frac{1}{\sqrt{n}}
I_\zs{n}(\psi_\zs{j})\,.
\end{equation}
\noindent
Through this representation  we introduce the following condition for the noise process \eqref{sec:In.1}
which will be used for the improved estimation below.

$\H_\zs{2}$) {\sl There exists $p_\zs{0}\ge 1$ and $d_\zs{0}\ge 1$ such that
for any $p\ge p_\zs{0}$ and $d\ge d_\zs{0}$ there exists a  $\sigma$ - field $\cG$
for which
 the random variables
 $(\xi_\zs{j})_\zs{1\le j\le d}$
 are the $\cG$-conditionally Gaussian with the covariance matrix
\begin{equation}\label{sec:Imp.6-1}
\G
=
\G_\zs{Q}
=\left(
\E_\zs{Q}\,\xi_\zs{i}\,\xi_\zs{j}|\cG)
\right)_\zs{1\le i,j\le d}
\end{equation}
and for some nonrandom
constant $l_\zs{*}=l_\zs{*}(p,d)>0$
\begin{equation}\label{sec:Imp.6-1-0}
\inf_\zs{Q\in\cQ_\zs{n}}\quad
\left(
\tr \,\G
-
\,
\lambda_\zs{max}(\G)
\right)
\geq l_\zs{*}
\quad \mbox{a.s.}\,,
\end{equation}
where $\tr\,\G$ is the trace and $\lambda_\zs{max}(\G)$ is the maximal eigenvalue of the matrix $\G$.
}

\noindent  Now we give an example
for this condition.

\begin{proposition}
\label{Pr.sec:Ex-CnD-H-2}
The conditional
covariance matrix
\eqref{sec:Imp.6-1}
constructed
for the  model \eqref{sec:In.1} -- \eqref{sec:Ex.1}
with respect to the $\cG=\sigma\{z_\zs{t}\,,\,t\ge 0\}$
on the basis of the trigonometric functions \eqref{sec:In.5}
for any odd $p\ge 3$
and the distribution family
\eqref{sec:Ex.01-1} -- \eqref{sec:Ex.01-2}
satisfies condition $\H_\zs{2})$
with $l_\zs{*}=\underline{\varrho}_\zs{n}(d-6)/2$
and
\begin{equation}
\label{d-0-Df-1}
d_\zs{0}=\inf\left\{d\ge 7\,:\,5+\ln d\le \check{a} d\right\}
\,,\qquad
\check{a}=\frac{1-e^{-a_\zs{max}}}{4a_\zs{max}}
\,.
\end{equation}
\end{proposition}

\noindent
Given condition $\H_\zs{2})$
 for the first $d$ Fourier coefficients in \eqref{sec:Imp.4},
 we will use
 the improved estimation method
proposed
for parametric models
 in
\cite{Pchelintsev2013}. We
define
\begin{equation}\label{sec:Imp.12}
\theta^{*}_\zs{j}=
\left(
1
-
g(j)
\right)
\,\wh{\theta}_\zs{j}
\quad\mbox{and}\quad
g(j)=\frac{\c_\zs{n}}{|\wh{\theta}|_\zs{d}}\, \Chi_\zs{\{1\le j\le d\}}
\,,
\end{equation}
where $\wh{\theta}=(\wh{\theta}_\zs{j})_\zs{1\le j\le p}\in\bbr^{p}$ and
 the norm $\vert x\vert^{2}_\zs{d}=\sum^{d}_\zs{j=1}\,x^{2}_\zs{j}$
for any vector $x=(x_\zs{j})_\zs{1\le j\le p}$ from $\bbr^{p}$ with $p\ge d$. Here
\begin{equation}
\label{coeff-n-c-n-10-13}
\c_\zs{n}=
\frac{l_\zs{*}}{n\left(r+\sqrt{\bar{d}\varkappa^{*}}\right)}\,,
\qquad
\bar{d}=\frac{d}{n}\,,
\end{equation}
and the coefficient $\varkappa^{*}$ is defined in condition $\H_\zs{1})$.
The positive parameter
 $r$  may be dependent of $n$, i.e. $r=r_\zs{n}$, and such that
\begin{equation}\label{sec:Imp.12+1}
\lim_\zs{n\to\infty}\,
n^{-\epsilon}\,r_\zs{n}
=0
\quad\mbox{for any}\quad
\epsilon>0\,.
\end{equation}

\noindent
To compare the estimators
\eqref{sec:Imp.3} and  \eqref{sec:Imp.12}
we put
\begin{equation}
\label{dif-risk-1-2}
\Delta_\zs{Q,S}:=\E_\zs{Q,S}\,
\vert\theta^{*}-\theta\vert^{2}_\zs{d}
-
\E_\zs{Q,S}\,
\vert\wh{\theta}-\theta\vert^{2}_\zs{d}\,,
\end{equation}
where $\theta^{*}=(\theta^{*}_\zs{j})_\zs{1\le j\le p}$ and $\wh{\theta}=(\wh{\theta}_\zs{j})_\zs{1\le j\le p}$ are the vectors in $\bbr^{p}$.
Now we study this difference.

\begin{theorem}\label{Th.sec:Imp.1}
Assume that the function $S$ is Lipschitzian. Then for any
$n\ge 1$, $d\ge 2$ and $p\ge 1$ for which conditions $\H_\zs{1})$ -- $\H_\zs{2})$ hold
\begin{equation}
\label{sec:Imp.11+1-n}
\Delta^{*}
=
\sup_{Q\in\cQ_\zs{n}}\,\sup_\zs{\Vert S\Vert\le r}
\Delta_\zs{Q,S}
\le
-\c^2_\zs{n}+\frac{2\sqrt{d}\phi_* L}{p}\c_\zs{n}
\,,
\end{equation}
where
$$
L=\sup_\zs{0\le s,t\le 1}
\frac{\vert S(t)-S(s)\vert}{\vert t-s\vert}
\,.
$$
\end{theorem}

\proof
First note that
$$
\E_\zs{Q,S} |\theta^{*}-\bar{\theta}|_\zs{d}^2=
\E _\zs{Q,S}|\wh{\theta}-\bar{\theta}|_\zs{d}^2+
\c^{2}_\zs{n}
-
2
\c_\zs{n}
\E_\zs{Q,S}\,\sum^{d}_\zs{j=1}J_\zs{j}
\,,
$$
where  $\bar{\theta}=(\bar{\theta}_\zs{j})_\zs{1\le j\le p}$,
$$
J_\zs{j}=\E_\zs{Q,S}\left(\iota_\zs{j}(\wh{\theta})(\wh{\theta}_\zs{j}-\bar{\theta}_\zs{j})|\cG\right)
\quad\mbox{and}\quad
\iota_\zs{j}(x)=\frac{x_\zs{j}}{|x|_\zs{d}}\,.
$$
In view of condition $\H_\zs{2})$ the vector $(\wh{\theta}_\zs{j})_\zs{1\le j\le d}$
is the $\cG$-conditionally Gaussian vector in $\bbr^{d}$  with mean  vector $\m=(\bar{\theta}_\zs{j})_\zs{1\le j\le d}$
and covariance matrix $n^{-1} \G$, then we obtain
\begin{equation*}
J_\zs{j}
=\int_{\mathbb{R}^d}\,\iota_\zs{j}(x)(x_\zs{j}-\bar{\theta}_\zs{j})\p(x|\cG)
\d x
\end{equation*}
with
\begin{equation*}
\p(x|\cG)=\frac{n^{d/2}}{(2\pi)^{d/2}\sqrt{\det(\G)}}
\exp\left(-\frac{n(x-\m)'\,\G^{-1}(x-\m)}{2}\right)
\,.
 \end{equation*}
Changing the variables by
$x=\s(u)=\G^{1/2}u + \m$
 and denoting by
$\v_\zs{ij}$ the $(i,j)$-th\ element of  $\G^{1/2}$,
we get
\begin{equation*}\label{sec:2.2c}
J_\zs{j}=\frac{n^{d/2}}{(2\pi)^{d/2}}
\sum_{l=1}^{d}
\v_\zs{j,l}\,
\int_{\bbr^{d}}\wt{\iota}_\zs{j}(u)u_{l}\exp\left(-\frac{n|u|^{2}_\zs{d}}{2}\right) \d u\,,
 \end{equation*}
where
$\wt{\iota}_\zs{j}(u)=\iota_\zs{j}(\s(u))$.
Furthermore,  the integrating by parts over the variable $u_\zs{l}$
implies that
$$
\frac{n^{d/2}}{(2\pi)^{d/2}}
\int_{\bbr^{d}}\wt{\iota}_\zs{j}(u)u_{l}\exp\left(-\frac{n|u|^{2}_\zs{d}}{2}\right) \d u
=\frac{1}{n}
\sum_\zs{k=1}^{d}
\,\v_\zs{k,l}\,\int_\zs{\bbr^{d}}\,H_\zs{j,k}(\s(u))\p_\zs{0}(u)\d u
$$
where $H_\zs{j,k}(z)=\partial \iota_\zs{j}(z)/\partial z_\zs{k}$,
$\p_\zs{0}$ is the Gaussian density in $\bbr^{d}$ with the parameters $(0,n^{-1}I_\zs{d})$ and $I_\zs{d}$ is the identity matrix of order $d$.
Taking into account that
$$
H_\zs{j,k}(z)=\frac{\vert z\vert^{2}_\zs{d}\Chi_\zs{\{k=j\}}-z_\zs{j}z_\zs{k}}{\vert z\vert^{3/2}_\zs{d}}
$$
and $\g_\zs{j,k}=\sum_{l=1}^{d}
\sum_\zs{k=1}^{d}\,
\v_\zs{j,l}
\,\v_\zs{k,l}$ is the $(j,k)$-th element of the matrix $\G$,
we find
$$
J_\zs{j}=
\E_\zs{Q,S}\left(
\frac{\g_\zs{j,j}}{\vert\wh{\theta}\vert_\zs{d}}
-
\frac{\sum^{d}_\zs{k=1}\g_\zs{j,k}\wh{\theta}_\zs{j}\wh{\theta}_\zs{k}}{\vert \wh{\theta}\vert^{3/2}_\zs{d}}
\vert
\cG
\right)\,.
$$
In view of the inequality
$\sum^{d}_\zs{k,j=1}\g_\zs{j,k}\wh{\theta}_\zs{j}\wh{\theta}_\zs{k}\le \lambda_{max}(\G) |\wh{\theta}|^{2}_\zs{d}$
and condition $\H_\zs{2})$
we have
\begin{align*}
\E_\zs{Q,S} |\theta^{*}-\bar{\theta}|_\zs{d}^2
&-\E_\zs{Q,S} |\wh{\theta}-\bar{\theta}|_\zs{d}^2 =
\c^{2}_\zs{n}
-
2
\c_\zs{n}
n^{-1}\E_\zs{Q,S}\,\left(\frac{\tr \G}{|\wh{\theta}|_\zs{d}}-
\frac{\sum^{d}_\zs{k,j=1}\g_\zs{j,k}\wh{\theta}_\zs{j}\wh{\theta}_\zs{k}}{|\wh{\theta}|_\zs{d}^3}\right)
\\[2mm]
&
\le \c^{2}_\zs{n}
-
\frac{2\c_\zs{n}\,l_\zs{*}}{n}
\,
\E_\zs{Q,S}\,\frac{1}{|\wh{\theta}|_\zs{d}}
\,.
\end{align*}
Besides, using the Jensen inequality, we can estimate the last expectation from below as
$$
\E_\zs{Q,S}\,\frac{1}{|\wh{\theta}|_\zs{d}}
=\E_\zs{Q,S}\,
\frac{1}{|\bar{\theta}+n^{-1/2}\xi|_\zs{d}}
\geq
\frac{1}{|\bar{\theta}|_\zs{d}+n^{-1/2}\E_\zs{Q,S}|\xi|_\zs{d}}
\,,
$$
where the vector  $\xi=(\xi_\zs{j})_\zs{1\le j\le p}$ and its components  are defined in \eqref{sec:Imp.4}.
From \eqref{sec:In.3}
we  get
$$
\E_\zs{Q,S}|\xi|^{2}_\zs{d}
=\sum^{d}_\zs{j=1}\,\E_\zs{Q,S}\,\xi^{2}_\zs{j}
\le \frac{\varkappa_\zs{Q}}{n}\sum^{d}_\zs{j=1}\int^{n}_\zs{0}\phi^{2}_\zs{j}\d t
\le \varkappa^*\,d\,.
$$
Since $\vert \bar{\theta}\vert_\zs{d}\le \Vert S\Vert$, then for $\Vert S\Vert\le r$
$$
\E_\zs{Q,S}\,|\wh{\theta}|^{-1}_\zs{d}\geq
\left(r+\sqrt{d\varkappa^*/n}\right)^{-1}
=\frac{1}{r+\sqrt{\bar{d}\varkappa^*}}
$$
and, therefore,
\begin{align}
\nonumber
\E_\zs{Q,S}\, |\theta^{*}-\bar{\theta}|_\zs{d}^2
\,&\le \E_\zs{Q,S} |\wh{\theta}-\bar{\theta}|_\zs{d}^2 + \c^{2}_\zs{n}
-
\frac{2\c_\zs{n}\,(d-1)l_\zs{*}}{n(r+\sqrt{\bar{d}\varkappa^*})}\\[2mm]\label{UpBnd.1}
&=
\E_\zs{Q,S} |\wh{\theta}-\bar{\theta}|_\zs{d}^2 - \c^{2}_\zs{n}
\,.
\end{align}
The applying here the representation \eqref{discoeff} implies
\begin{align*}
\E_\zs{Q,S}\, |\theta^{*}-\bar{\theta}|_\zs{d}^2&=
\E_\zs{Q,S}\,|\theta^{*}-\theta|_\zs{d}^2
- 2\sum^{d}_\zs{j=1}\,\E_\zs{Q,S} (\theta^{*}_\zs{j}-\theta_\zs{j}) h_\zs{j}+|h|^2_\zs{d}
\\[2mm]
&\le \E_\zs{Q,S} |\wh{\theta}-\theta|_\zs{d}^2 - \c^{2}_\zs{n}
- 2\sum^{d}_\zs{j=1}\,\E_\zs{Q,S} (\wh{\theta}_\zs{j}-\theta_\zs{j}) h_\zs{j}+|h|^2_\zs{d}
\,.
\end{align*}
So, the risks difference
\eqref{dif-risk-1-2} can be estimated from above as
\begin{equation}\label{sec:Prf.2}
\Delta_\zs{Q,S}
\leq -\c^2_n+2
\sum^{d}_\zs{j=1}\,\E_\zs{Q,S} (\theta^{*}_\zs{j}-\wh{\theta}_\zs{j}) h_\zs{j}
\,.
\end{equation}
Taking into account that for any $\varepsilon>0$
\begin{equation*}
2\vert ab\vert\leq \varepsilon a^2 +\varepsilon^{-1}b^2\,,
\end{equation*}
we obtain
$$
2
\sum^{d}_\zs{j=1}\,\E_\zs{Q,S} (\theta^{*}_\zs{j}-\wh{\theta}_\zs{j}) h_\zs{j}
\leq \varepsilon \E_\zs{Q,S} |\theta^*-\wh{\theta}|_d^2+ \varepsilon^{-1} |h|_d^2
= \varepsilon \c^2_n+\varepsilon^{-1} |h|_d^2\,.
$$
Now, from definition of the coefficients $h_\zs{j}$ in  \eqref{discoeff} and by the Cauchy-Bunyakovsky-Schwarz inequality, we have the estimate
\begin{align*}
|h|_d^2\,&
=\sum_{j=1}^d \left(\sum_{k=1}^p \int_{t_{k-1}}^{t_{k}}\phi_\zs{j}(t_{k})\left(S(t)-S(t_{k})\right)dt \right)^2 \\[3mm]
&\,\leq \sum_{j=1}^d \sum_{k=1}^p \int_{t_{k-1}}^{t_{k}}\phi_\zs{j}^2(t_{k})\d t\,
\sum_{k=1}^p \int_{t_{k-1}}^{t_{k}}\left(S(t)-S(t_{k})\right)^2 \d t \le \frac{d\phi_*^2 L^2}{p^2}
\,.
\end{align*}
Using this in
\eqref{sec:Prf.2}, we get
$$
\Delta_\zs{Q,S}
\leq -(1-\varepsilon)\c^2_n+
\varepsilon^{-1}
\frac{d\phi_*^2 L^2}{p^2}
\,.
$$
Minimizing this upper bound over $\varepsilon>0$, we come to the inequality \eqref{sec:Imp.11+1-n}. Hence
Theorem \ref{Th.sec:Imp.1}.
\endproof

\noindent
Let now
\begin{equation}
\label{pp00-LB}
p_\zs{0}=\frac{2\sqrt{d}\phi_\zs{*} L}{\c_\zs{n}}
=
\frac{2\phi_\zs{*} L(r+\sqrt{\bar{d}\varkappa^{*}})n\sqrt{d}}{l_\zs{*}}\,,
\end{equation}
where $L$ is the Lipschitz constant for $S$ defined in \eqref{sec:Imp.11+1-n}.

\begin{remark}\label{Re;sec:Imp.1}
Note that uniformly over $p$  for $p> p_\zs{0}$ the upper bound in \eqref{sec:Imp.11+1-n} is negative, i.e. $\Delta^{*}<0$.
This  means that non-asymptotically, i.e. for the bounded number observations $n$,
 the estimate \eqref{sec:Imp.12}   outperforms in mean square accuracy the estimate \eqref{sec:Imp.3}.
\end{remark}

\begin{remark}\label{Re;sec:Imp.2}
For the model \eqref{sec:In.1} -- \eqref{sec:Ex.1}
with $\underline{\varrho}_\zs{n}\equiv \underline{\varrho}_\zs{*}>0$
asymptotically as $n,d\to\infty$
for
\begin{equation}
\label{cnd-freq-1}
p\ge 2p_\zs{0}
\approx
\frac{n}{\sqrt{d}}
\end{equation}
 we obtain that
$$
\Delta^{*}\le
-
\frac{\underline{\varrho}^{2}_\zs{*}(d-6)^{2}}{2\left(r+\sqrt{\bar{d}\varkappa^{*}}\right)^{2}\,n^{2}}
\approx -\left(\frac{d}{n}\right)^{2}
\,.
$$
We recall that in the parametric case the upper bound in \eqref{sec:Imp.11+1-n}
asymptotically goes to zero as $n^{-2}$ (see, \cite{KonevPergamenshchikovPchelintsev2014, Pchelintsev2013}) since the parametric dimension is fixed, but in this case
the dimension $d\to\infty$. Therefore, Theorem \ref{Th.sec:Imp.1}  implies that for the nonparametric models
with incomplete observations
 the improvement effect is more significant than in the parametric models uniformly over
observation frequency.
\end{remark}

\section{Improved model selection}\label{sec:Mo}

In this section we construct  model selection procedures
based on the improved estimators \eqref{sec:Imp.12}. To this end, by the same way as
 in \cite{KonevPergamenshchikov2015} we define a class of weighted least squares estimates
 $(S^{*}_\zs{\gamma})_\zs{\gamma\in\Gamma}$ as
\begin{equation}\label{sec:Imp.5}
\wh{S}_\zs{\gamma}(t)=\sum^{p}_\zs{j=1}\gamma(j)\wh{\theta}_\zs{j}\psi_\zs{j}(t)\,,
\end{equation}
where $\wh{\theta}_\zs{j}$ is defined in \eqref{sec:Imp.4}, the weights $\gamma=(\gamma(j))_{1\leq j\leq p}\in\bbr^{p}$ belong to some
finite set $\Gamma\subset [0,1]^{p}$. For this set we denote
\begin{equation}\label{sec:Imp.6}
\nu= \#(\Gamma)
\quad\mbox{and}\quad
\nu_\zs{*}=\max_{\gamma\in \Gamma} \, \sum_{j=1}^p\gamma(j)\,,
\end{equation}
where $\#(\Gamma)$ is the number of the vectors $\gamma$ in $\Gamma$.
In the sequel we assume that all vectors from $\Gamma$ satisfies the following condition.

\medskip

$\H_\zs{3}$) {\sl  Assume that the set $\Gamma$ is such  for any vector $\gamma\in\Gamma$ there exists
 some fixed integer $d=d(\gamma)\le p$
 such that their first $d$ components
are equal to one, i.e. $\gamma(j)=1$ for $1\le j\le d$ for any $\gamma\in\Gamma$.  Moreover, we assume that
the parameters $\nu$ and $\nu_\zs{*}$ are functions of $n$, i.e. $\nu=\nu(n)$ and $\nu_\zs{*}=\nu_\zs{*}(n)$,
such that  for any  $\epsilon>0$
\begin{equation}
\label{Gamma-Cnds-1}
\lim_\zs{n\to\infty}
\frac{\nu}{n^{\epsilon}}= 0
\quad\mbox{and}\quad
\lim_\zs{n\to\infty}
\frac{\nu_\zs{*}}{n^{1/3+\epsilon}}= 0
\,.
\end{equation}
}

\noindent
Now we define shrinkage estimates for $S$
\begin{equation}\label{sec:Imp.11}
S^{*}_\zs{\gamma}(t)=\sum^{p}_\zs{j=1}\gamma(j)
\theta^{*}_\zs{j}\psi_\zs{j}(t)\,,
\end{equation}
where the shrinkage estimators  $\theta^{*}_\zs{j}$ are defined in  \eqref{sec:Imp.12}
with the parameter $d$ given in condition $\H_\zs{3}$). In this case the threshold $\c_\zs{n}$
is the function of $\gamma$ also, i.e. $\c_\zs{n}=\c_\zs{n}(\gamma)$.  In the sequel we will use the maximum of this
sequence, i.e.
\begin{equation}\label{Unif-c-1}
\c^{*}_\zs{n}= n\,\max_\zs{\gamma\in\Gamma}\,\c^{2}_\zs{n}(\gamma)
\,.
\end{equation}

\noindent
$\H_\zs{4}$) {\sl
Assume that $\lim_\zs{n\to\infty}n^{-\epsilon} \c^{*}_\zs{n}=0$ for any $\epsilon>0$.}

\noindent
Now
we compare the estimators \eqref{sec:Imp.5} and \eqref{sec:Imp.11}.

\begin{theorem}\label{Th.sec:Imp.1-1-n}
Assume that the function $S$ is Lipschitzian. Then for any
$n\ge 3$, $d\ge 2$, $p\ge  1$ and $\gamma\in\Gamma$
 for which the conditions $\H_\zs{1})$ --  $\H_\zs{3})$ hold
\begin{equation}
\label{sec:Imp.11+1}
\sup_\zs{Q\in\cQ_\zs{n}}\,\sup_\zs{\Vert S\Vert\le r}
\left(
\cR_\zs{Q}(S^{*}_\zs{\gamma},S)-\cR_\zs{Q}(\wh{S}_\zs{\gamma},S)
\right)
\le-\c^{2}_\zs{n}
\,,
\end{equation}
where the coefficient $\c_\zs{n}$ is defined in
\eqref{coeff-n-c-n-10-13}.
\end{theorem}

\proof  First, introducing the projection
\begin{equation}
\label{prj-S-121}
\bar{S}=\sum^{p}_\zs{j=1}\bar{\theta}_\zs{j}\psi_\zs{j}(t)\,,
\end{equation}
we have
$$
\Vert
S^{*}_\zs{\gamma}
-
S
\Vert^2
=
\Vert
S^{*}_\zs{\gamma}
-
\bar{S}
\Vert^{2}
+\Vert S-\bar{S}\Vert^{2}
\quad\mbox{and}\quad
\Vert
\wh{S}_\zs{\gamma}
-
S
\Vert^2
=
\Vert
\wh{S}_\zs{\gamma}
-
\bar{S}
\Vert^{2}
+\Vert S-\bar{S}\Vert^{2}
\,.
$$
Therefore, by condition $\H_\zs{3})$, we get
$$
\cR_\zs{Q}(S^{*}_\zs{\gamma},S)-\cR_\zs{Q}(\wh{S}_\zs{\gamma},S)
=
\E_\zs{Q,S}\, |\theta^{*}-\bar{\theta}|_\zs{d}^2- \E_\zs{Q,S} |\wh{\theta}-\bar{\theta}|_\zs{d}^2\,.
$$
Using the inequality \eqref{UpBnd.1}, we obtain the bound \eqref{sec:Imp.11+1}.
Hence Theorem \ref{Th.sec:Imp.1-1-n}.
\endproof

\begin{remark}\label{Re.sec:MdSel.1-001}
We would like to emphasize that the accuracy improvement provided by the upper bound \eqref{sec:Imp.11+1}   for
the weighted least squares estimate \eqref{sec:Imp.11}
is uniform over the observation frequency $p\ge 1$.  Indeed, we obtained the same upper bound as for the
estimation problem on the complete observations (see, Theorem 3.1 in \cite{PchelintsevPergamenshchikov2018b}).
This is the new
 effect in the improvement nonparametric estimation theory.
\end{remark}

\noindent
To construct  model selection procedures
we need to impose some stability conditions
 introduced in \cite{PchelintsevPergamenshchikov2018a}
 for the noise sequence in
\eqref{sec:Imp.4}.

\noindent $\C_\zs{1})$ {\it There exists a proxy variance
$\sigma_\zs{Q}> 0$
such that for any $\epsilon>0$
\begin{equation*}
\label{sec:C1EQ--00}
\lim_\zs{n\to\infty}\frac{1}{n^{\epsilon}}\,
\sup_\zs{Q\in\cQ_\zs{n}}\,
\sum^{p}_\zs{j=1}\,\left|
\E_\zs{Q}\,\xi^{2}_\zs{j}
-
\sigma_\zs{Q}
\right|=0
 \,.
\end{equation*}
}
\noindent
Further we will need the following function
\begin{equation}
\label{sec:C1EQ}
\L_\zs{1}(Q)=\sum^{p}_\zs{j=1}\,
\left|
\E_\zs{Q}\,\xi^{2}_\zs{j}
-
\sigma_\zs{Q}
\right|
 \,.
\end{equation}

\noindent
To study the total noise variance in \eqref{sec:Imp.4} we denote
\begin{equation}
\label{sec:C2EQ}
\L_{2}(Q)=\sup_\zs{x\in \cB_\zs{p}}
\E_\zs{Q}\,
\left(
\sum^{p}_\zs{j=1}\,
x_\zs{j}\,\wt{\xi}_\zs{j}
\right)^{2}
\quad\mbox{and}\quad
\wt{\xi}_\zs{j}
=
\xi^{2}_\zs{j}-\E_\zs{Q} \xi^{2}_\zs{j}
\,,
\end{equation}
where $\cB_\zs{p}=\{x\in\bbr^{p}\,:\,\vert x\vert_\zs{p}\le 1\}$ and
 $\vert \cdot\vert_\zs{p}$ is Euclidean norm  in $\bbr^{p}$.

\noindent $\C_\zs{2})$ {\it Assume that  the sequence
$\L^{*}_\zs{2,n}=\sup_\zs{Q\in\cQ_\zs{n}}\L_{2}(Q)$ is such that
\begin{equation*}
\label{sec:C2EQ--001}
\lim_\zs{n\to\infty}\frac{\L^{*}_\zs{2,n}}{n^{\epsilon}}=0
\quad\mbox{for any}\quad
\epsilon>0
\,.
\end{equation*}
}

\noindent
As is shown in Propositions \ref{Pr.sec:L1-10} - \ref{Pr.sec:L2-11}, the
trigonometric basis \eqref{sec:In.5} provides the both conditions $\C_\zs{1})$
and $\C_\zs{2})$ for  the
noise process \eqref{sec:Ex.1}
with the properties \eqref{sec:Ex.01-1} -- \eqref{sec:Ex.01-2}.

To choose a weight vector $\gamma$ in $\Gamma$ we use the  quadratic estimation accuracy
$$
\Er\,(\gamma)=\|S^*_\zs{\gamma}-S\|^2.
$$
By the definition
\eqref{sec:Imp.11} and the projection \eqref{prj-S-121} we can represent this accuracy as
\begin{align}\nonumber
\Er(\gamma)&=\sum^{p}_\zs{j=1}(\gamma(j)\theta^*_\zs{j}-\bar{\theta}_\zs{j})^{2}
+\Vert S-\bar{S}\Vert^{2}
\\[2mm] \label{sec:Mo.1}
&=\sum^{p}_\zs{j=1}\,\gamma^2(j)(\theta^*_\zs{j})^2\,-
2\,\sum^{p}_\zs{j=1}\,\gamma(j)\theta^*_\zs{j}\,\bar{\theta}_\zs{j}\,+\,
\|S\|^2\,.
\end{align}
It is clear, that  to choose a good weighted  vector
 one needs to minimize this function over $\gamma$
 in $\Gamma$. The coefficients
 $(\bar{\theta}_\zs{j})_\zs{1\le j\le p}$
 are unknown, and we need to replace in this function the terms
$\theta^*_\zs{j}\,\bar{\theta}_\zs{j}$ by their estimators.
Similar to
\cite{KonevPergamenshchikov2015}
in this case we use the estimators defined as
\begin{equation}\label{sec:Mo.2}
\wt{\theta}_\zs{j}=
\theta^{*}_\zs{j}\,\wh{\theta}_\zs{j}-\frac{\wh{\sigma}_\zs{n}}{n}\,,
\end{equation}
where $\wh{\sigma}_\zs{n}$ is an estimate
for the proxy variance $\sigma_\zs{Q}$ defined in the condition $\C_\zs{1})$.
As in \cite{KonevPergamenshchikov2015}
we define it
through the estimators of the Fourier coefficients on the
trigonometric basis \eqref{sec:In.5}, i.e.
\begin{equation}\label{sec:Mo.3}
\wh{\sigma}_\zs{n}=
\frac{n}{p}
\sum_\zs{j=[\sqrt{n}]+1}^n \wh{\t}^2_\zs{j}
\quad\mbox{and}\quad
\wh{\t}_\zs{j}=
\frac{1}{n}\sum^{np}_\zs{l=1}\,
\Trg_\zs{j}(t_\zs{l})\,(y_\zs{t_\zs{l}}-y_\zs{t_\zs{l-1}})\,.
\end{equation}
\noindent
According to the model selection approach for the dependent observations
 (see, for example, in \cite{GaltchoukPergamenshchikov2009a, KonevPergamenshchikov2015}),
 when we modify the quadratic accuracy
 one needs to penalize this modification by adding in \eqref{sec:Mo.1}  a special positive penalty term.
 In this case the penalization means that  the additional positive term makes the minimization problem  more difficult. So, for any $\gamma\in\Gamma$
  we set the objective function as
 \begin{equation}\label{sec:Mo.4}
J(\gamma)\,=\,\sum^{p}_\zs{j=1}\,\gamma^2(j)(\theta^{*}_\zs{j})^2\,-
2\,\sum^{p}_\zs{j=1}\,\gamma(j)\,\wt{\theta}_\zs{j}\,
+\,\rho\,\wh{P}_\zs{n}(\gamma)\,,
\end{equation}
where $\rho$ is some positive penalization level,
$\wh{P}_\zs{n}(\gamma)$ is the penalty term defined as
\begin{equation}\label{sec:Mo.5}
\wh{P}_\zs{n}(\gamma)=\frac{\wh{\sigma}_\zs{n}\,|\gamma|^2_\zs{p}}{n}
\,.
\end{equation}
 In the case, when the value of $\sigma_\zs{Q}$ in $\C_\zs{1})$ is known, we have
$\wh{\sigma}_\zs{n}=\sigma_\zs{Q}$ and
\begin{equation}\label{sec:Mo.9}
P_\zs{n}(\gamma)=\frac{\sigma_\zs{Q}\,|\gamma|^2_\zs{p}}{n}\,.
\end{equation}

\noindent
Finally, we define the improved model selection procedure as
\begin{equation}\label{sec:Mo.7}
\gamma^*=\mbox{argmin}_\zs{\gamma\in\Gamma}\,J(\gamma)
\quad\mbox{and}\quad
S^*=S^*_\zs{\gamma^*}\,.
\end{equation}
It will be noted that $\gamma^*$ exists because
 $\Gamma$ is a finite set. If the
minimizing sequence in \eqref{sec:Mo.7} $\gamma^*$ is not
unique, one can take any minimizer.

\noindent $\C_\zs{3}$) {\it Assume that for any $n\ge 1$
\begin{equation}
\label{sigma*-nn}
\varsigma^{*}=\varsigma^{*}_\zs{n}=\sup_\zs{Q\in\cQ_\zs{n}}
\sigma_\zs{Q}
<\infty
\quad\mbox{and}\quad
\varsigma_\zs{*}=\varsigma_\zs{*,n}
=
\inf_\zs{Q\in\cQ_\zs{n}}
\sigma_\zs{Q}>0
\end{equation}
and for any  $\epsilon>0$ the ratio $\varsigma^{*}_\zs{n}/n^{\epsilon}\to 0$ as $n\to\infty$.
}

\noindent
We use this condition to construct the special set $\Gamma$ of weight vectors $(\gamma(j))_\zs{j\ge 1}$ as it is proposed  in \cite{GaltchoukPergamenshchikov2009a, GaltchoukPergamenshchikov2009b}
 for which we will study the asymptotic properties of the model selection procedure \eqref{sec:Mo.7}.
 For this we consider the following grid
\begin{equation*}\label{sec:Imp.7}
\cA_\zs{n}=\{1,\ldots,k^*\}\times\{r_1,\ldots,r_m\}\,,
\end{equation*}
where  $r_i=i\varepsilon$, $i=\overline{1,\,m}$ with $m=[1/\varepsilon^2]$. We assume that the
parameters $k^*\ge 1$ and $0<\varepsilon\le 1$ are
functions of $n$, i.e. $ k^*=k^*(n)$ and
$\varepsilon=\varepsilon(n)$, such that
\begin{equation}\label{sec:Ga.1-305_2019}
\lim_\zs{n\to \infty}\,
\left(
\frac{1}{k^{*}(n)}+
\frac{k^{*}(n)}{\ln n}
\right)
=0
\quad\mbox{and}\quad
\lim_\zs{n\to \infty}\,
\left(\varepsilon(n)
+
\frac{1}{n^{\b}\varepsilon(n)}
\right)
=0
\end{equation}
for any $\b>0$. One can take, for example, for $0<\ve<1$
\begin{equation}\label{sec:Ga.1-00}
\varepsilon(n)=1/\ln (n+1)
\quad\mbox{and}\quad
k^{*}(n)=k^{*}_\zs{0}+\sqrt{\ln (n+1)}\,,
\end{equation}
where $k^{*}_\zs{0}\ge 0$ is some fixed constant.
 \noindent
For each $\alpha=(\beta,r)\in\cA_\zs{n}$ we introduce the weight
sequence $\gamma_\zs{\alpha}=(\gamma_\zs{\alpha}(j))_\zs{j\ge 1}$
as
\begin{equation}\label{sec:Imp.9}
\gamma_\zs{\alpha}(j)=\Chi_\zs{\{1\le j\le j_\zs{*}(\alpha)\}}+
\left(1-(j/\omega_\alpha)^\beta\right)\, \Chi_\zs{\{ j_\zs{*}(\alpha)<j\le
\omega_\alpha\}}\,,
\end{equation}
where $j_\zs{*}(\alpha)=\omega_\zs{\alpha}/\ln (n+1)$,
$$
\omega_\zs{\alpha}=\left(\frac{(\beta+1)(2\beta+1)}{\pi^{2\beta}\beta}\,r\,v_\zs{n}\right)^{1/(2\beta+1)}
\quad\mbox{and}\quad
v_\zs{n}=n/\varsigma^{*}
\,.
$$

\noindent Finally, we set
\begin{equation}\label{sec:Imp.10_Lambda}
\Gamma\,=\,\{\gamma_\zs{\alpha}\,,\,\alpha\in\cA_\zs{n}\}\,.
\end{equation}

\noindent As we will see later, the model selection procedure  \eqref{sec:Mo.7}
with this weight set will  be  efficient in adaptive setting.

\begin{remark}\label{Re.sec:MdSel.1-0}

As it is shown in \cite{KonevPergamenshchikov2015}
 the weight coefficients \eqref{sec:Imp.10_Lambda}   satisfy condition $\H_\zs{3})$
with  $d=[j_\zs{*}(\alpha)]$. Moreover, for the  model \eqref{sec:In.1} -- \eqref{sec:Ex.1} the using of the trigonometric basis \eqref{sec:In.5} and
Proposition \ref{Pr.sec:Ex-CnD-H-2} provide condition  $\H_\zs{4})$ with the family
\eqref{sec:Imp.10_Lambda}.
\end{remark}

\begin{remark}\label{Re.sec:MdSel.1}
As it is shown in \cite{Pinsker1981} in the case when the regularity parameters
of the function $S$ are known,
 the weighted least squares estimate
\eqref{sec:Imp.5}  with the weights of the form \eqref{sec:Imp.9}  is asymptotically efficient. In this  paper
we use these weights to provide through Theorem \ref{Th.sec:Imp.1-1-n} the efficient property for the improved model selection procedure
\eqref{sec:Mo.7} in the adaptive case, i.e. when the regularity properties are unknown.
\end{remark}

\section{Main results}\label{sec:MainRs}

In this Section we obtain
the sharp oracle inequalities for the quadratic risk \eqref{sec:In.4} and robust risk \eqref{sec:In.6} of proposed procedure.
Then, on the basis of these inequalities, we will establish the robust efficiency property in adaptive setting.

\subsection{Oracle inequalities}\label{sec:Orineq}

\noindent
First, we obtain the oracle inequalities for the risks \eqref{sec:In.4}.
\begin{theorem}\label{sec:Mo.Th.1}
Assume that conditions $\C_\zs{1})$ - $\C_\zs{3})$ hold. Then, for any $n\geq 3$, $p\ge 1$
and the weight family $\Gamma$ for which conditions $\H_\zs{1})$ - $\H_\zs{4})$ hold also
and for any $0<\rho<1/2$
\begin{equation}\label{Or-INEQ-12-13}
\cR_\zs{Q}(S^{*},S)\,\le\, \frac{1+5\rho}{1-\rho}
\min_\zs{\gamma\in\Gamma} \cR_\zs{Q}(S^*_\zs{\gamma},S)
+
\frac{\U_\zs{n}
\left(
1+\nu_\zs{*}\E_\zs{Q,S}|\wh{\sigma}_\zs{n}-\sigma_\zs{Q}|
\right)
}{\rho n}
\,,
\end{equation}
where
 the coefficient $\U_\zs{n}$ is such that
\begin{equation}
\label{termB_rest}
\lim_\zs{n\to\infty}\frac{\U_\zs{n}}{n^{\epsilon}}=0
\quad\mbox{for any}\quad \epsilon>0
\,.
\end{equation}
\end{theorem}

\noindent
Now we study the estimate \eqref{sec:Mo.3}. We need the following condition for the observation
frequency.

\vspace{2mm}

$\D)${\sl Assume that the observation frequency $p$ is a function of $n$, i.e. $p=p(n)$ such that
$1\le p\le n$ and
$\lim_\zs{n\to\infty} n^{\epsilon-5/6}p=\infty$ for any $\epsilon>0$.
}

\begin{proposition}
\label{Pr.sigma-est-1}
Assume that the conditions
 $\C_\zs{1})$ -- $\C_\zs{3})$ hold for the trigonometric basis \eqref{sec:In.5}.
 Moreover, assume also that the unknown function $S$ has the square integrated derivative $\dot{S}$ and the conditions $\H_\zs{1})$ and $\D)$ hold.
 Then for $n\ge 3$ and $\sqrt{n}<p\le n$
\begin{equation}
\label{UpBnd-sigma-101-1}
\E_\zs{Q,S}|\wh{\sigma}_\zs{n}-\sigma_\zs{Q}\vert\le \K_\zs{n} n^{-1/3}
\left(1+
\|\dot{S}\|^2
\right)\,,
\end{equation}
where the term $\K_\zs{n}>0$ is such that 
$\lim_\zs{n\to\infty}\,n^{-\epsilon}\,\K_\zs{n}
=0$ for any $\epsilon>0$ .
\end{proposition}

\noindent
Theorem \ref{sec:Mo.Th.1} and Proposition \ref{Pr.sigma-est-1} imply the oracle inequalities for the robust risks \eqref{sec:In.6}.
\begin{theorem}\label{Th.sec:2.3}
Assume that conditions $\C_\zs{1})$ - $\C_\zs{3})$, $\D)$ hold  and the function $S$ has the square integrable derivative $\dot{S}$. Then, for any $n\geq 3$, $p\ge 1$
and the weight family $\Gamma$ for which the conditions $\H_\zs{1})$ - $\H_\zs{4})$ hold also
and for any $0<\rho<1/2$\begin{align*}\label{sec:Mo.20}
\cR^{*}_\zs{n}(S^{*},S)\,\le\,
\frac{1+5\rho}{1-\rho} \min_\zs{\gamma\in\Gamma}
\cR_n^{*}(S^*_\zs{\gamma},S) +\frac{1}{\rho n}\,\U_\zs{n}(1+\|\dot{S}\|^2)\,,
\end{align*}
where the term $\U_\zs{n}$ satisfies the property
\eqref{termB_rest}.
\end{theorem}

\begin{remark}\label{Re.sec:MdSel.3344-01}
Conditions
$\C_\zs{1})$ -- $\C_\zs{3})$ are usually required  to obtain nonasymptotic sharp oracle inequalities for the statistical models
with the dependent observations
through the methods developed in
\cite{GaltchoukPergamenshchikov2009a, KonevPergamenshchikov2009a}.
Condition $\D)$ is basic to provide  the property \eqref{UpBnd-sigma-101-1}.
Moreover, conditions $\H_\zs{1})$ - $\H_\zs{4})$ are used to study  the properties of the shrinkage estimators
\eqref{sec:Imp.11}.
\end{remark}

\begin{remark}\label{Re.sec:MdSel.2-34-1}
Main difference between the weight estimators \eqref{sec:Imp.5} and \eqref{sec:Imp.11}
is that generally the shrinkage estimators $\theta^{*}_\zs{j}$ depend on the weight vector $\gamma$. By this reason to show the inequality \eqref{Or-INEQ-12-13}
 we can't use directly
the method developed in \cite{KonevPergamenshchikov2015} for the weighted least squares estimates \eqref{sec:Imp.5}. For such estimator one needs to develop
a new special analytical tool.
\end{remark}

\subsection{Asymptotic efficiency}\label{sec:Ae}

In order to study the asymptotic efficiency we define the following functional Sobolev ball
\begin{equation*}\label{sec:Ae.1}
W_\zs{k,\r}=\left\{f\in\cC^{(k)}_\zs{per}[0,1]\,:\,
\sum_\zs{j=0}^{k}\,\|f^{(j)}\|^2\le \r\right\}\,,
 \end{equation*}
where $\r>0$ and $k\ge 1$ are
some unknown parameters, $\cC^{k}_\zs{per}[0,1]$ is the space of
 $k$ times differentiable $1$-periodic $\bbr\to\bbr$ functions
 such that  for any $0\le i \le k-1$ the periodic boundary conditions are satisfied, i.e. $f^{(i)}(0)=f^{(i)}(1)$.
To analyze  the efficiency properties for the model selection procedure \eqref{sec:Mo.7}
one needs to compare this estimator with all possible estimation methods. To this end
we denote by $\Xi_\zs{n}$ the set of all estimators $\wh{S}_\zs{n}$ based on the observations $(y_\zs{t})_\zs{0\le t\le n}$
i.e. any
$\sigma\{y_\zs{t}\,,\,0\le t\le n\}$ measurable functions.
To obtain the lower bound for the risks  we will use  the process  \eqref{sec:In.1}
with $\xi_\zs{t}=\varsigma^{*}w_\zs{t}$, i.e. the white noise model with the intensity  \eqref{sigma*-nn}.
We denote its distribution by $Q^{*}$.

\begin{theorem}\label{Th.sec: Ae.1}
Assume that $Q^{*}\in\cQ_\zs{n}$.
Then robust risks \eqref{sec:In.6} are bounded from below in the following sense
 \begin{equation}\label{sec:LB.0-5-2}
\liminf_\zs{n\to\infty}\,
v_\zs{n}^{2k/(2k+1)}
\inf_\zs{\wh{S}_\zs{n}\in\Xi_\zs{n}}
\,
\sup_\zs{S\in W_\zs{k,\r}}\,\cR^{*}_\zs{n}(\wh{S}_\zs{n},S)
\,
\ge l_\zs{k}(\r)
\,,
 \end{equation}
where $v_\zs{n}=n/\varsigma^{*}$ and  $l_\zs{k}(\r)\,=\,((2k+1)\r)^{1/(2k+1)}\,
\left(k/(\pi (k+1))\right)^{2k/(2k+1)}$.
\end{theorem}

\noindent
We show that this lower bound is sharp in the following sense. To provide the efficient property for  the model selection procedure
\eqref{sec:Mo.7} we choose the parameter $\rho$ as a function of $n$, i.e. $\rho=\rho(n)$ such that
\begin{equation}
\label{rho-1Cnd}
\lim_\zs{n\to\infty}\,\rho_\zs{n}=0
\quad\mbox{and for any }\, \epsilon>0\qquad
\lim_\zs{n\to\infty}\,n^{\epsilon}\rho_\zs{n}=\infty\,.
\end{equation}

\begin{theorem}\label{Th.sec: Ae.2}
Assume that conditions $\C_\zs{1})$ - $\C_\zs{3})$, $\D)$, $\H_\zs{1})$ - $\H_\zs{4})$ hold.
Then robust risk of  the model selection procedure
 \eqref{sec:Mo.7} constructed on the trigonometric basis functions \eqref{sec:In.5}
  with the weight coefficients \eqref{sec:Imp.10_Lambda} and the penalty threshold satisfying the
 conditions \eqref{rho-1Cnd} is bounded from above as
 \begin{equation*}\label{sec:Ae.4}
\limsup_\zs{n\to\infty}\,v_\zs{n}^{2k/(2k+1)}
\sup_\zs{S\in W_\zs{k,\r}}\,\cR^{*}_\zs{n}(S^{*},S)
\,
\le l_\zs{k}(\r)
\,.
 \end{equation*}
\end{theorem}

\noindent
It is clear that these theorems
imply the following efficient property.
\begin{theorem}\label{Th.sec: Ae.1-EFF-1}
Assume that the conditions of Theorems
\ref{Th.sec: Ae.1} and  \ref{Th.sec: Ae.2} hold.
Then the model selection procedure
 \eqref{sec:Mo.7} constructed on the trigonometric basis functions \eqref{sec:In.5} with the weight coefficients \eqref{sec:Imp.10_Lambda}
 and the penalty threshold satisfying the
 conditions \eqref{rho-1Cnd}
is asymptotically efficient, i.e.
\begin{equation*}\label{sec:Ae.5--10}
\lim_\zs{n\to\infty}\,v_\zs{n}^{2k/(2k+1)}\,
\sup_\zs{S\in W_\zs{k,\r}}\,\cR^{*}_\zs{n}(S^{*},S)\,
= l_\zs{k}(\r)
\end{equation*}
and
$$
\lim_\zs{n\to\infty}\frac{\inf_\zs{\wh{S}_\zs{n}\in\Sigma_\zs{n}}
\,
\sup_\zs{S\in W_\zs{k,\r}}\,\cR^{*}_\zs{n}(\wh{S}_\zs{n},S)}{\sup_\zs{S\in W_\zs{k,\r}}\,\cR^{*}_\zs{n}(S^{*},S)}
=1\,.
$$
\end{theorem}
\noindent
Theorem \ref{Th.sec: Ae.1} is shown by the same way as Theorem 1
in \cite{KonevPergamenshchikov2009b}.
Theorem \ref{Th.sec: Ae.2} follows from Theorems \ref{Th.sec:2.3}, \ref{Th.sec:Imp.1-1-n} and Theorem 5.2
in \cite{KonevPergamenshchikov2015}.

\begin{remark}\label{Re.sec:Pinsk-1201}
The level $l_\zs{k}(\r)$ determining the lower bound \eqref{sec:LB.0-5-2} is the well-known Pinsker constant, obtained in
\cite{Pinsker1981} for the filtration problem in the model \eqref{sec:In.1} with the standard white noise $(\xi_\zs{t})_\zs{t\ge 0}$ that modeled by
the standard Brownian motion.  For the models with general semimartingale noise the lower bound is the same as for
the white noise model, but generally  the convergence rate is not the same.  In this case the convergence rate is given by
$\left(n/\varsigma^{*}_\zs{n}\right)^{-2k/(2k+1)}$ while in classical white noise model the convergence rate is $\left(n\right)^{-2k/(2k+1)}$.
So, if the upper variance threshold $\varsigma^{*}_\zs{n}$ tends to zero, the convergence rate is better than the classical one; if it tends to infinity, it is worse
and, if it is a constant, the rate is the same.
\end{remark}

\subsection{Statistical analysis for the big data model \eqref{sec:nparam-md-1}}
Now we apply our results for the high dimensional model \eqref{sec:nparam-md-1}. We assume that the functions $(\u_\zs{j})_\zs{1\le j\le q}$ are
orthonormal in $\L_\zs{2}[0,1]$, i.e. 
\begin{equation}
\label{bg-dt-1}
S(t)=\sum^{q}_\zs{j=1}\beta_\zs{j}\u_\zs{j}(t)
\,.
\end{equation}
We use the estimator \eqref{sec:Imp.11}
to estimate the parameters $\beta=(\beta_\zs{j})_\zs{1\le j\le q}$ as
$$
\beta^{*}_\zs{\gamma}=(\beta^{*}_\zs{\gamma,j})_\zs{1\le j\le q}
\quad\mbox{and}\quad
\beta^{*}_\zs{\gamma,j}=(\u_\zs{j},S^{*}_\zs{\gamma})
\,.
$$
Moreover, we use the selection model procedure \eqref{sec:Mo.7}
as
\begin{equation}
\label{selmd-1}
\beta^{*}=(\beta^{*}_\zs{j})_\zs{1\le j\le q}
\quad\mbox{and}\quad
\beta^{*}_\zs{j}=(\u_\zs{j},S^{*})
\,.
\end{equation}
It is clear that
$$
\vert \beta^{*}_\zs{\gamma}-\beta \vert^{2}_\zs{q}=
\sum^{q}_\zs{j=1}(\beta^{*}_\zs{\gamma,j}-\beta_\zs{j})^2
=\Vert S^{*}_\zs{\gamma}-S\Vert^{2}
\quad\mbox{and}\quad
\vert \beta^{*}-\beta \vert^{2}_\zs{q}
=\Vert S^{*}-S\Vert^{2}\,.
$$
Therefore, Theorem \ref{Th.sec:2.3} implies
\begin{theorem}\label{Th.sec:bgd-1.1}
Assume that conditions $\C_\zs{1})$ - $\C_\zs{3})$, $\D)$ hold  and the function \eqref{bg-dt-1} has the square integrable derivative $\dot{S}$. Then, for any $n\geq 3$, $p\ge 1$
and the weight family $\Gamma$ for which conditions $\H_\zs{1})$ - $\H_\zs{4})$ hold also
and for any $0<\rho<1/2$
\begin{align}\nonumber
\sup_\zs{Q\in\cQ_\zs{n}}\E_\zs{Q,\beta} \vert \beta^{*}-\beta \vert^{2}_\zs{q}
&\le\,
\frac{1+5\rho}{1-\rho} \min_\zs{\gamma\in\Gamma}
\,\sup_\zs{Q\in\cQ_\zs{n}}\E_\zs{Q,\beta} \vert \beta^{*}-\beta \vert^{2}_\zs{q}
\\[3mm]\label{sec:Bgd.1-1}
&
 +\frac{1}{\rho n}\,\U_\zs{n}(1+\|\dot{S}\|^2)\,,
\end{align}
where the term $\U_\zs{n}$ satisfies the property
\eqref{termB_rest}.
\end{theorem}

\noindent
Theorems
\ref{Th.sec: Ae.1} and  \ref{Th.sec: Ae.2}
imply the efficiency property for the estimate \eqref{selmd-1} based on
the model selection procedure
 \eqref{sec:Mo.7}
constructed on the trigonometric basis functions \eqref{sec:In.5} with the weight coefficients \eqref{sec:Imp.10_Lambda} and
the penalty threshold satisfying the
 conditions \eqref{rho-1Cnd}.
\begin{theorem}\label{Th.sec: bgd121.1-EFF}
Assume that the conditions of Theorems
\ref{Th.sec: Ae.1} and  \ref{Th.sec: Ae.2} hold.
Then the estimate \eqref{selmd-1}
is asymptotically efficient, i.e.
\begin{equation*}\label{sec:Ae.5--10}
\lim_\zs{n\to\infty}\,v_\zs{n}^{2k/(2k+1)}\,
\sup_\zs{S\in W_\zs{k,\r}}\,
\sup_\zs{Q\in\cQ_\zs{n}}\E_\zs{Q,\beta} \vert \beta^{*}-\beta \vert^{2}_\zs{q}
= l_\zs{k}(\r)
\end{equation*}
and
$$
\lim_\zs{n\to\infty}\frac{\inf_\zs{\wh{\beta}_\zs{n}\in\Xi_\zs{n}}
\,
\sup_\zs{S\in W_\zs{k,\r}}\,\sup_\zs{Q\in\cQ_\zs{n}}\E_\zs{Q,\beta} \vert \wh{\beta}_\zs{n}-\beta \vert^{2}_\zs{q}}{\sup_\zs{S\in W_\zs{k,\r}}\,\sup_\zs{Q\in\cQ_\zs{n}}\E_\zs{Q,\beta} \vert \beta^{*}-\beta \vert^{2}_\zs{q}}
=1\,,
$$
where  $\Xi_\zs{n}$ is the set of all possible estimators for the vector $\beta$.
\end{theorem}

\begin{remark}\label{Re.sec:BgD-1.1}
In the estimator \eqref{sec:Mo.7}  doesn't use the dimension $q$ in \eqref{bg-dt-1}. Moreover, it can be equal to $+\infty$.
In this case it is impossible to use neither LASSO method nor Danzig selector.
\end{remark}

\begin{remark}\label{Re.sec:BgD-1.2}
If in addition the functions $(\u_\zs{j})_\zs{1\le j\le q}$ are orthonormal
on the grid $\{t_\zs{1},\ldots t_\zs{p}\}$, i.e. possess the  property \eqref{innerprod}, as, for example, in the trigonometric basis case
\eqref{sec:In.5},  we can use these functions
in \eqref{sec:Imp.3}.
\end{remark}

\section{Monte-Carlo simulations}\label{sec:Sim}

In this section we give the results of numerical simulations to assess the performance and improvement of
the proposed model selection procedure \eqref{sec:Mo.7}.
We simulate the model \eqref{sec:In.1} with
$1$-periodic functions $S$ of the forms
\begin{equation}\label{sec:Sim_Sign_11}
S_1(t)=t\,\sin(2\pi t)+t^2(1-t)\cos(4\pi t)
\end{equation}
and
\begin{equation}\label{sec:Sim_Sign_11-1}
S_2(t)=\sum^{+\infty}_\zs{j=1}\frac{1}{1+j^{3}}\sin(2\pi j t)
\end{equation}
on $[0,\,1]$ and the Ornstein -- Uhlenbeck -- L\'evy noise process $\xi_\zs{t}$ is defined as
$$
\d \xi_\zs{t}=-\xi_\zs{t} \d t + 0.5\,\d w_\zs{t}+0.5\, \d z_\zs{t}\,,\quad
z_\zs{t}=\sum^{N_\zs{t}}_\zs{j=1}\,Y_\zs{j}
\,,
$$
where $N_\zs{t}$ is  a homogeneous  Poisson process of intensity $\lambda=1$ and
$(Y_\zs{j})_\zs{j\ge1}$ is i.i.d.
$\cN(0,\,1)$ sequence (see, for example, \cite{KonevPergamenshchikov2015}).

We use the model selection procedure   \eqref{sec:Mo.7} with the weights \eqref{sec:Imp.9} in which
$k^*=100+\sqrt{\ln (n+1)}$, $r_i=i/\ln (n+1)$, $m=[\ln^2 (n+1)]$, $\varsigma^*=0.5$ and $\rho=(3+\ln n)^{-2}$.
We define the empirical risk as
$$
\cR(S^*,\,S)=\frac{1}{p}\sum_\zs{j=1}^p
\wh{\E}
\Delta^{2}_\zs{n}(t_\zs{j})
\quad\mbox{and}\quad
 \wh{\E}
 \Delta^{2}_\zs{n}(t)
 =
\frac{1}{N}\sum_\zs{l=1}^N
\Delta^{2}_\zs{n,l}(t)
\,,
$$
where
$\Delta_\zs{n}(t)=S_\zs{n}^*(t)-S(t)$
and
$\Delta_\zs{n,l}(t)=S_\zs{n,l}^*(t)-S(t)$ is the deviation for the $l$-th replication. In this example we take $p=10001$ and
 $N = 1000$.


\begin{table}[ht]
\label{Tab1}
\caption{The sample quadratic risks for different optimal $\gamma$}
\begin{center}
\begin{tabular}{lrrrr}
  \hline
   $n$              & 100 & 200 & 500 & 1000 \\ \hline
  $\cR(S^*_\zs{\gamma^*},\,S_1)$    & 0.0819 & 0.0319 & 0.0098 & 0.0051 \\ 
  $\cR(\wh{S}_\zs{\wh{\gamma}},\,S_1)$ & 0.0787 & 0.0479 & 0.0287 & 0.0178 \\ 
  $\cR(\wh{S}_\zs{\wh{\gamma}},\,S_1)/\cR(S^*_\zs{\gamma^*},\,S_1)$ & 0.9 & 1.5 & 2.9 & 3.5 \\
   $\cR(S^*_\zs{\gamma^*},\,S_2)$    & 1.2604 & 0.6979 & 0.2571 & 0.0269 \\ 
  $\cR(\wh{S}_\zs{\wh{\gamma}},\,S_2)$ & 3.8215 & 2.2983 & 0.9728 & 0.1398 \\ 
  $\cR(\wh{S}_\zs{\wh{\gamma}},\,S_2)/\cR(S^*_\zs{\gamma^*},\,S_2)$ & 3.1 & 3.3 & 3.8 & 5.2 \\
  \hline
\end{tabular}
\end{center}
\end{table}

\begin{table}[ht]
\label{Tab2}
\caption{The sample quadratic risks for the same optimal $\wh{\gamma}$}
\begin{center}
\begin{tabular}{lrrrr}
  \hline
   $n$              & 100 & 200 & 500 & 1000 \\ \hline
  $\cR(S^*_\zs{\wh{\gamma}},\,S_1)$    & 0.0671 & 0.0397 & 0.0195 & 0.0097 \\ 
  $\cR(\wh{S}_\zs{\wh{\gamma}},\,S_1)$ & 0.0787 & 0.0479 & 0.0287 & 0.0178 \\ 
  $\cR(\wh{S}_\zs{\wh{\gamma}},\,S_1)/\cR(S^*_\zs{\wh{\gamma}},\,S_1)$ & 1.2 & 1.2 & 1.5 & 1.8 \\
$\cR(S^*_\zs{\wh{\gamma}},\,S_2)$    & 2.7992 & 1.6650 & 0.4322 & 0.0531 \\ 
  $\cR(\wh{S}_\zs{\wh{\gamma}},\,S_2)$ & 3.8215 & 2.2983 & 0.9728 & 0.1398\\ 
  $\cR(\wh{S}_\zs{\wh{\gamma}},\,S_2)/\cR(S^*_\zs{\wh{\gamma}},\,S_2)$ & 1.4 & 2.2 & 2.3 & 2.6 \\
  \hline
\end{tabular}
\end{center}
\end{table}

Table 1 gives the values for the sample risks of the improved estimate \eqref{sec:Mo.7}
and the model selection procedure based on the weighted LSE (3.15) from \cite{KonevPergamenshchikov2012} for different numbers
of observation period $n$. Table 2 gives the values for the sample risks of the the model selection procedure based on the weighted LSE (3.15) from \cite{KonevPergamenshchikov2012} and it's improved version for different numbers
of observation period $n$.


\begin{figure}[h!]
\label{fig1}
\begin{minipage}[h]{0.49\linewidth}
\center{\includegraphics[width=0.9\textwidth]{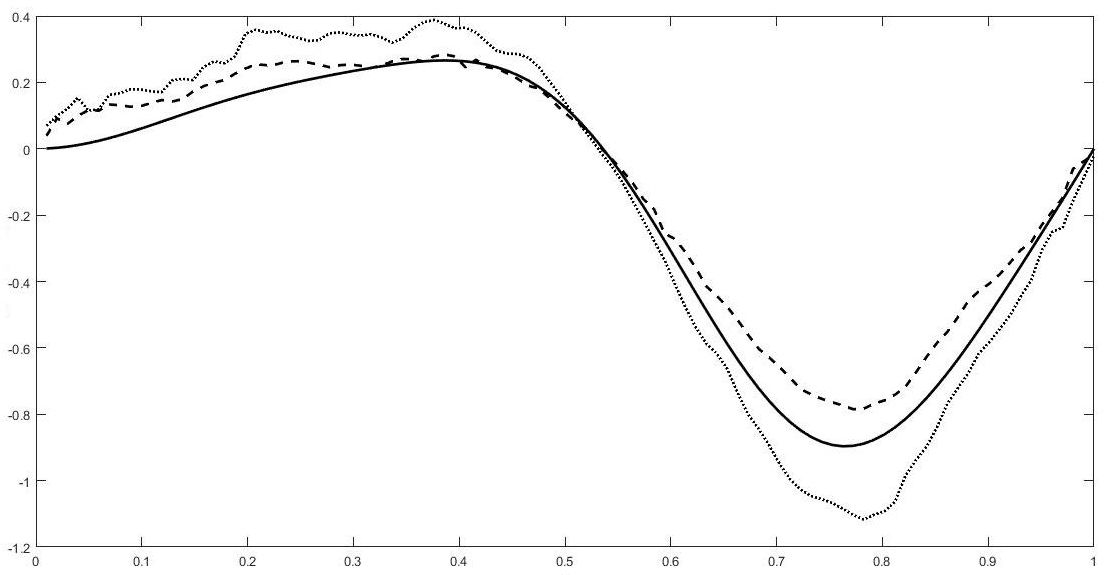} \\ a)}
\end{minipage}
\hfill
\begin{minipage}[h]{0.49\linewidth}
\center{\includegraphics[width=0.7\textwidth]{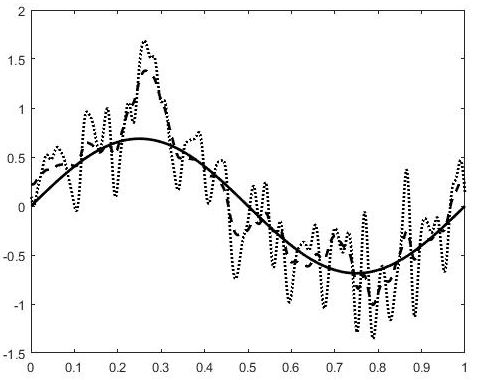} \\ b)}
\end{minipage}
\caption{Behavior of the regression functions and their estimates for $n=100$ (\, a) -- for the function $S_1$ and b) -- for the function $S_2$).}
\end{figure}

\begin{figure}[h!]
\label{fig1-23-1}
\begin{minipage}[h]{0.49\linewidth}
\center{\includegraphics[width=0.9\textwidth]{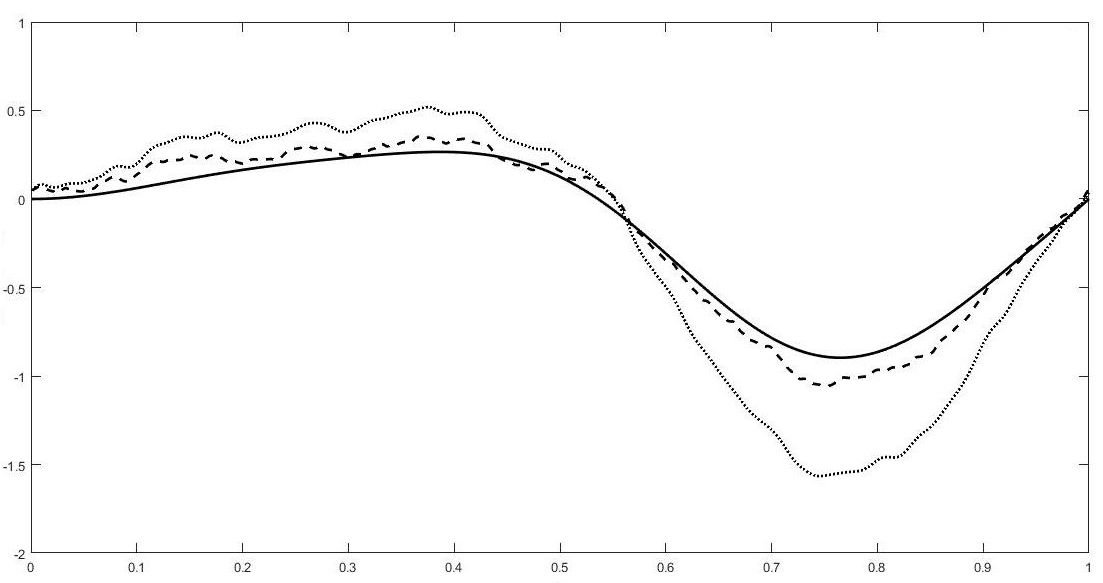} \\ a)}
\end{minipage}
\hfill
\begin{minipage}[h]{0.49\linewidth}
\center{\includegraphics[width=0.7\textwidth]{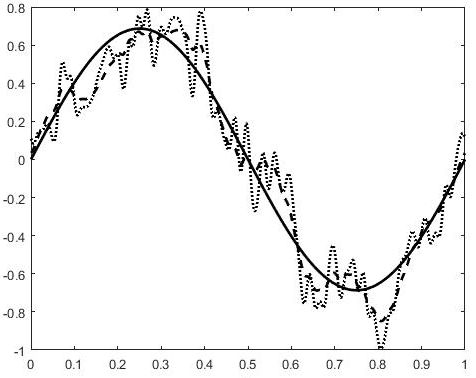} \\ b)}
\end{minipage}
\caption{Behavior of the regressions function and their estimates for $n=500$ (\, a) -- for the function $S_1$ and b) -- for the function $S_2$).}
\end{figure}

\begin{figure}[h!]
\label{fig1-23-2-2}
\begin{minipage}[h]{0.49\linewidth}
\center{\includegraphics[width=0.9\textwidth]{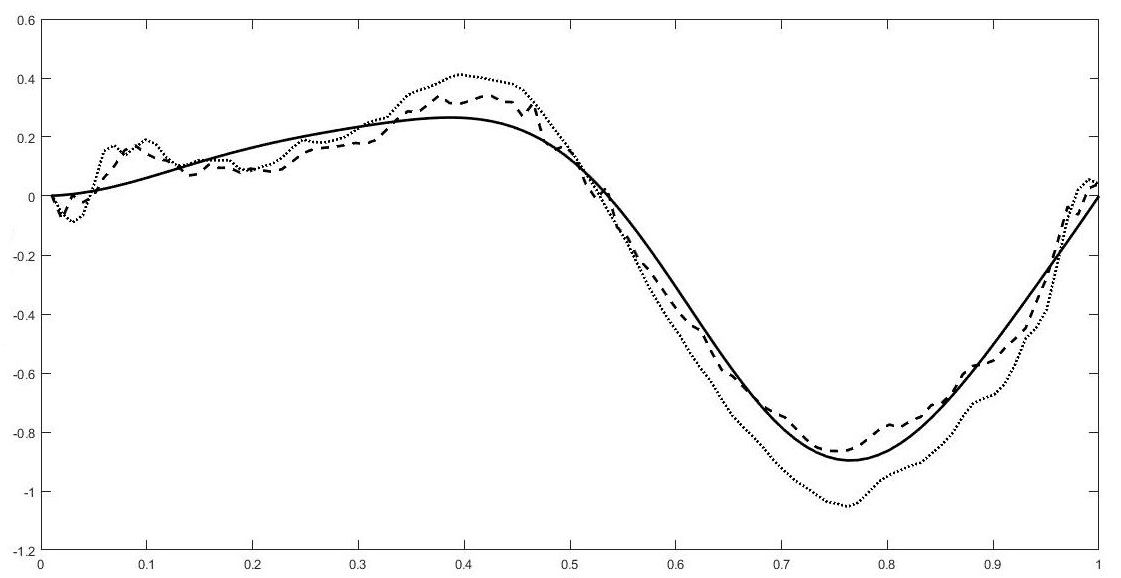} \\ a)}
\end{minipage}
\hfill
\begin{minipage}[h]{0.49\linewidth}
\center{\includegraphics[width=0.7\textwidth]{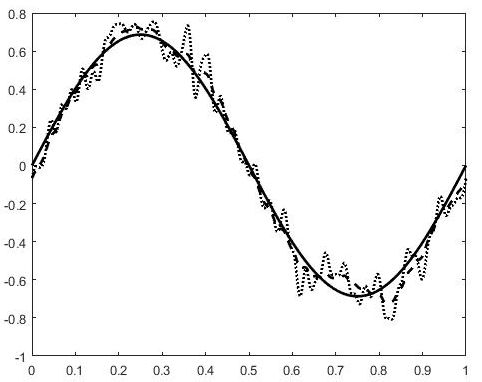} \\ b)}
\end{minipage}
\caption{Behavior of the regression functions and their estimates for $n=1000$ (\, a) -- for the function $S_1$ and b) -- for the function $S_2$).}
\end{figure}

%
%
%
%

\begin{remark}
Figures show the behavior of
the procedures  \eqref{sec:Imp.5} and \eqref{sec:Mo.7}   depending
on the values of observation periods $n$.
The bold lines are the functions \eqref{sec:Sim_Sign_11} and \eqref{sec:Sim_Sign_11-1},
the continuous  lines are the  model selection
procedures based on the
least squares estimates $\wh{S}$
 and the dashed lines are the improved model selection procedures $S^*$.
 From the Table 2 for the same $\gamma$ with various observations
numbers $n$ we can conclude that
theoretical result on the improvement effect \eqref{sec:Imp.11+1} is confirmed
by the numerical simulations.
Moreover, for the proposed shrinkage procedure, from the Table 1 and Figures 1--3,
we can conclude that the benefit is considerable for non large  $n$.
\end{remark}

\section{Properties of the model \eqref{sec:nparam-md-1} -\eqref{sec:Ex.1}}
\label{sec:Main-prps}

In this section
 we study the  noise process \eqref{sec:Ex.1}
with the conditions
\eqref{sec:Ex.1-00_mPi}-\eqref{sec:Ex.01-2}.  More precisely,
we check conditions $\C_\zs{1})$ and $\C_\zs{2})$ for  the trigonometric basis
\eqref{sec:In.5}.
First we study condition $\C_\zs{1})$.

\begin{proposition}\label{Pr.sec:L1-10}
There exists a constant $\l>0$ such that for any $n\ge 1$
\begin{equation}
\label{UbnD-1-0-2}
\sup_\zs{Q\in\cQ_\zs{n}}\L_{1}(Q)
\le  \,\l\,
\varsigma^*\,,
\end{equation}
where the parameter $\varsigma^{*}$ is defined in \eqref{sec:Ex.01-1}.
 \end{proposition}

\proof
First
we note that Proposition \ref{Pr.sec:Stc.1} implies
\begin{equation}\label{sec:Pr.14-1}
\E_\zs{Q,S} \xi^{2}_\zs{j}= \sigma_\zs{Q} \left( 1+\tau_\zs{j}\right)\,,
\end{equation}
where $\tau_\zs{j}=n^{-1}\,\int^{n}_\zs{0}\,\psi_\zs{j}(t)
\check{\varepsilon}_\zs{j}(t)\,\d t$
and
$\check{\varepsilon}_\zs{j}(t)=a
\int^{t}_\zs{0}\,e^{a(t-s)}\,
\psi_\zs{j}(s)\,(1+e^{2as})
\d s$.
It is easy to check
that $|\tau_\zs{1}|\le 2$.
For $j\ge 2$ we represent $\tau_\zs{j}$ and $\check{\varepsilon}_\zs{j}(t)$ as
$$
\tau_\zs{j}=\tau_\zs{1,j}+\tau_\zs{2,j}
\quad\mbox{and}\quad
\check{\varepsilon}_\zs{j}(t)=
\varepsilon_\zs{1,j}(t)+
\varepsilon_\zs{2,j}(t)
\,,
$$
where
$
\varepsilon_\zs{1,j}(t)=
a
\int^{t}_\zs{0}\,e^{a(t-s)}\,
\psi_\zs{j}(s)\,
\d s$,
$\varepsilon_\zs{2,j}(t)=
a
\int^{t}_\zs{0}\,e^{a(t+s)}\,
\psi_\zs{j}(s)\,
\d s$

\begin{equation} \label{sec:Pr.14-2}
\tau_\zs{1,j}=\frac{1}{n}
\int^{n}_\zs{0}\,\psi_\zs{j}(t)
\varepsilon_\zs{1,j}(t)\,\d t
\quad\mbox{and}\quad
\tau_\zs{2,j}=\frac{1}{n}\int^{n}_\zs{0}\,\psi_\zs{j}(t)
\varepsilon_\zs{2,j}(t)\,\d t\,.
\end{equation}
In view of $\psi_{j}(t)=\sum_{k=1}^{np}\Trg_\zs{j}(t_k)\Chi_\zs{(t_{k-1},t_k]}(t)$ the values $\tau_\zs{1,j}$ can be represented as
$$
\tau_\zs{1,j}
=\frac{a}{n}\,
\sum^{np}_\zs{k=1}\,
\sum^{k}_\zs{l=1}\,
\Trg_\zs{j}(t_\zs{k})\,
\Trg_\zs{j}(t_\zs{l})\,
\varkappa_\zs{k,l}\,,
$$
where
$$
\varkappa_\zs{k,l}=\int^{t_\zs{l}}_\zs{t_\zs{l-1}}\,
\left(
\int^{t_\zs{k}}_\zs{t_\zs{k-1}}\,
(
e^{a(t-s)}
\,\Chi_\zs{\{s\le t\}}
)
\d t
\right)
\,\d s\,.
$$
Integration yields
$\varkappa_\zs{k,k}=\varkappa^{*}_\zs{0}(p)=
(e^{a/p}-1-a/p)a^{-2}$
and for $1\le l<k$ the coefficients
$\varkappa_\zs{k,l}=\varkappa^{*}_\zs{1}(p)\,e^{a t_\zs{k-l}}$
with
$\varkappa^{*}_\zs{1}(p)=a^{-2}\,
(e^{-a/p}-1)\,(1-e^{a/p})$. Therefore,
\begin{equation}\label{sec:Pr.15}
\tau_\zs{1,j}
=a\varkappa^{*}_\zs{0}(p)
p +a\,
\varkappa^{*}_\zs{1}(p)\,p^{2}\,A_\zs{j,p}
\,,
\end{equation}
where
$$
A_\zs{j,p}=
\frac{1}{np^{2}}\,\sum^{np-1}_\zs{l=1}\,
\upsilon_\zs{j,l}\,,
\quad
\upsilon_\zs{j,l}=\sum^{np}_\zs{k=l+1}\,
\Trg_\zs{j}(t_\zs{k})\,
\Trg_\zs{j}(t_\zs{l})\,
e^{at_\zs{k-l}}
\,.
$$
It should be noted that
\begin{equation}\label{sec:Pr.15-00}
0\le \varkappa^{*}_\zs{0}(p)
\le \frac{1}{2p^{2}}
\quad\mbox{and}\quad
0\le \varkappa^{*}_\zs{1}(p)
\le \,\frac{e^{a_\zs{max}}}{p^{2}}
\,.
\end{equation}
From the definition of the functions $(\Trg_\zs{j})_\zs{2\le j\le p}$
it follows that
$$
\upsilon_\zs{j,l}=\sum^{np-l}_\zs{k=1}\,
\cos(\alpha_\zs{j} t_\zs{k})\,
e^{a t_\zs{k}}
+
(-1)^{j}
\sum^{np-l}_\zs{k=1}\,
\cos(\alpha_\zs{j}\,t_\zs{k+2l})\,
e^{at_\zs{k}}
:=\wt{\upsilon}_\zs{j,l}
+
(-1)^{j}
\wh{\upsilon}_\zs{j,l}
\,,
$$
where $\alpha_\zs{j}=2\pi[j/2]$. So, setting
$$
\wt{A}_\zs{j,p}=\frac{1}{np^{2}}\,\sum^{np-1}_\zs{l=1}\,
\wt{\upsilon}_\zs{j,l}
\quad\mbox{and}\quad
\wh{A}_\zs{j,p}=\frac{1}{np^{2}}\,\sum^{np-1}_\zs{l=1}\,
\wh{\upsilon}_\zs{j,l}\,,
$$
we can represent $A_\zs{j,p}$ in the form
 $A_\zs{j,p}=
\wt{A}_\zs{j,p}+(-1)^{j}
\wh{A}_\zs{j,p}$, in which
\begin{equation}\label{sec:Pr.15-0}
\wt{A}_\zs{j,p}=\,\frac{1}{p}
\cR e
\frac{\wt{\q}_\zs{j}}{1-\wt{\q}_\zs{j}}
\,
-\wt{B}_\zs{j,p}\quad\mbox{and}\quad
\wh{A}_\zs{j,p}=
\frac{1-e^{an}}{np^{2}}
\,
\wh{B}_\zs{j}
\,.
\end{equation}
Here
$$
\wt{B}_\zs{j,p}=\,
\frac{1-e^{a n}}{np^{2}}
\,
\cR e
\frac{\wt{\q}_\zs{j}}{(1-\wt{\q}_\zs{j})^{2}}
\quad\mbox{and}\quad
\wh{B}_\zs{j}=\,
\cR e\,
\frac{\wt{\q}_\zs{j}\wh{\q}_\zs{j}}
{(1-\wt{\q}_\zs{j})(1-\wh{\q}_\zs{j})}\,,
$$
where
$\wt{\q}_\zs{j}=e^{a/p+i\alpha_\zs{j}/p}$ and
$\wh{\q}_\zs{j}=e^{-a/p+i\alpha_\zs{j}/p}$.
 First, note that
$$
\cR e
\frac{1}{1-\wt{\q}_\zs{j}}
=
\frac{1-e^{a/p}\cos(\alpha_\zs{j}/p)}
{1-2e^{a/p}\cos(\alpha_\zs{j}/p)+e^{2a/p}}
=\frac{1-e^{a/p}+2e^{a/p}\sin^{2}(\alpha_\zs{j}/(2p))}
{(1-e^{a/p})^{2}+4e^{a/p}\sin^{2}(\alpha_\zs{j}a/(2p))}
\,.
$$
Taking into account that
$$
\sin(x)\ge 2x/\pi,\quad 0\le x\le\pi/2\,,
$$
 and
$$
\alpha_\zs{j}\ge 2\pi j/3\,,
\quad 2\le j\le p\,,
$$
we obtain for $2\le j\le p$
\begin{equation}\label{sec:Pr.15-0-1}
2\,\sin(\alpha_\zs{j}/(2p))
\ge j/p
\,.
\end{equation}
Therefore,
\begin{equation}\label{sec:Pr.15-2}
\frac{1}{p}\,
\left|\cR e
\frac{1}{1-\wt{\q}_\zs{j}}\right|
\,
\le
\frac{|a|/p}{4pe^{a/p}\sin^{2}(\alpha_\zs{j}/(2p))}
+\frac{1}{2p}
\,
\le
\frac{a_\zs{max} e^{a_\zs{max}}}{j^{2}}
\,+\,\frac{1}{2p}
\,.
\end{equation}
Now we estimate $\wt{B}_\zs{j,p}$. We have
\begin{align*}
\wt{B}_\zs{j,p}&=\frac{(1-e^{an})e^{a/p}\left(\cos(\alpha_\zs{j}/p)-
2e^{a/p}+e^{2a/p}\cos(\alpha_\zs{j}/p)\right)}
{np^{2}\left(1-2e^{a/p}\cos(\alpha_\zs{j}/p)+ e^{2a/p}\right)^{2}}\\[4mm]
&=
\frac{(1-e^{an})e^{a/p}
\left(
\left(1-e^{a/p}\right)^{2}
-
2\left(1+e^{2a/p}\right)
\sin^{2}(\alpha_\zs{j}/(2p))
\right)}{np^{2}\left(
\left(1-e^{a/p}\right)^{2}
+
4e^{a/p}\,\sin^{2}(\alpha_\zs{j}/(2p))
\right)^{2}}\,.
\end{align*}
Applying the inequality \eqref{sec:Pr.15-0-1} yields
\begin{align}\nonumber
|\wt{B}_\zs{j,p}|&\le
\frac{|a|^{2}}{16p^{4}e^{a/p}\sin^{4}(\alpha_\zs{j}/(2p))}
+
\frac{1}{4p^{2}e^{a/p}\sin^{2}(\alpha_\zs{j}/(2p))}
\\[2mm]\label{sec:Pr.15-3}
&
\le
\,e^{a_\zs{max}}\left(
a^{2}_\zs{max}\,j^{-4}
+
  j^{-2}
\right)\,.
\end{align}

\noindent By making use of the estimates \eqref{sec:Pr.15-2} and \eqref{sec:Pr.15-3}
in \eqref{sec:Pr.15-0} we obtain

\begin{equation}\label{sec:Pr.15-4}
\sup_\zs{n\ge 1}\,
|\wt{A}_\zs{j,p}|\le
\frac{1}{2 p}
+
\frac{(a_\zs{max}+1)e^{a^{max}}}{j^{2}}
+
\frac{a^{2}_\zs{max}e^{a^{max}}}{j^{4}}
\,.
\end{equation}
To estimate $\wh{A}_\zs{j,p}$ we represent $\wh{B}_\zs{j}$ as
$$
\wh{B}_\zs{j}=
\frac{\cos(\alpha_\zs{j}/p)}{2(\cos(\alpha_\zs{j}/p)-\ch(a/p))}\,,
$$
where $\ch(x)=(e^{x}+e^{-x})/2$. From here and
\eqref{sec:Pr.15-0-1}, it follows that
\begin{equation}\label{sec:Pr.15-5}
 |\wh{A}_\zs{j,p}|\le
\,
\le
\frac{1}{2np^{2}}
\left(
\frac{1}{1-\cos(\varpi_\zs{j}\Delta)}
\right)
=
\frac{1}{4np^{2}\sin^{2}(\varpi_\zs{j}\Delta/2)}
\le\,
\frac{1}{nj^{2}}
\,.
\end{equation}
Combining \eqref{sec:Pr.15-4} and \eqref{sec:Pr.15-5} yields
$$  p^{-1}
\sup_\zs{n\ge 1}\,
|A_\zs{j,p}|\le
\A_\zs{*}
\left(
  p^{-1}
+  j^{-2}
\right)
\,,\quad
\A_\zs{*}=1+(1+a_\zs{max}+a^{2}_\zs{max})e^{a_\zs{max}}\,.
$$
This and \eqref{sec:Pr.15}, in view of \eqref{sec:Pr.15-00}, implies that
\begin{equation}\label{sec:Pr.15-6}  j^{-2}
 \sup_\zs{n\ge 1}\,
|\tau_\zs{1,j}|\le
\,
\tau^{*}_\zs{1}\,
\left(
  p^{-1}
+
  j^{-2}
\right)\,,
\quad \tau^{*}_\zs{1}
=a_\zs{max}
\left(1
+
e^{a_\zs{max}}\,\A_\zs{*}
\right)
\,.
\end{equation}
It remains to estimate $\tau_\zs{2,j}$ in \eqref{sec:Pr.14-2}. First we note
that
$$
\tau_\zs{2,j}=\frac{a}{2 n}\,\iota^{2}_\zs{j,p}
\quad\mbox{and}\quad
\iota_\zs{j,p}=
\int^{n}_\zs{0}\,\psi_\zs{j}(t)\,e^{a t}\d t\,.
$$
It is easy to check  that for  $2\le j\le p$
\begin{equation}\label{sec:Pr.15-1}
\iota_\zs{j,p}=\frac{\sqrt{2}(e^{-a/p}-1)}{|a|}
\,\left(
1-e^{a n}
\right)
\,
\left(
\Upsilon_\zs{j}
\left(\frac{1}{1-\wt{\q}_\zs{j}}
\right)
-
1
\right)\,,
\end{equation}
where $\Upsilon_\zs{j}(z)=\cR e(z)$ for even $j$ and
 $\Upsilon_\zs{j}(z)=\cI m(z)$ for odd $j$.
For even $j$, in view of \eqref{sec:Pr.15-2}, one gets the inequality
\begin{equation}\label{sec:Pr.15-7}
|\iota_\zs{j,p}|\,
\le\,
\iota_\zs{*}\,
\left(
  p^{-1}
+
  j^{-2}
\right)\,,
\quad
\iota_\zs{*}=\sqrt{2}\,a_\zs{max}\,e^{2a_\zs{max}}+3\,e^{a_\zs{max}}\,.
\end{equation}
For odd $3\le j\le p$ one has the estimate
$$
\left|
\cI m \frac{1}{1-\wt{\q}_\zs{j}}
\right|
\le
\frac{p}{ j}\,,
$$
which implies that
\begin{equation}\label{sec:Pr.15-8}
|\iota_\zs{j,p}|\,
\le\,
\sqrt{2}\,e^{a_\zs{max}}
\,
\left(
  p^{-1}
+
  j^{-1}
\right)
\le
\,\iota_\zs{*}\,
\left(
  p^{-1}
+
  j^{-1}
\right)
\,.
\end{equation}
\noindent
Therefore,
$$
\sup_\zs{n\ge 1}\,
|\tau_\zs{2,j}|\le a_\zs{max}\,\iota^{2}_\zs{*}\,
\left(
p^{-2}
+
  j^{-2}
\right)\,,\quad 2\le j\le p\,.
$$
From here and the definition of $\tau_\zs{j}$ in \eqref{sec:Pr.14-1},
we obtain
$$
\sup_\zs{n\ge 1}\,
\sum^{p}_\zs{j=1}
\,|\tau_\zs{j}|
\le
\tau_\zs{*}
\left(
1
+
\sum_\zs{j\ge 1}\,
  j^{-2}
\right)\le 3\tau_\zs{*}\,.
$$
So, in view of \eqref{sec:Pr.14-1},
for any any $n\ge 1$ and $Q\in\cQ_\zs{n}$
the term  $\L_\zs{1}(Q)\le 3\sigma_\zs{Q}\tau_\zs{*}$ and $\tau_\zs{*}=a_\zs{max}\left(1+
e^{a_\zs{max}}+2\,(\sqrt{2} a_\zs{max} e^{a_\zs{max}}
+3)^{2}\,e^{2a_\zs{max}}\right)$.
Taking into account  that $\sigma_\zs{Q}\le \varsigma^*$, we get
the upper bound
\eqref{UbnD-1-0-2}. Hence Proposition~\ref{Pr.sec:L1-10}.
\endproof

\begin{proposition}\label{Pr.sec:L2-11}
There exists a constant $\l>0$ such that
for any  $n\ge 1$
\begin{equation}
\label{Up-Bnd-21_1}
\sup_\zs{Q\in\cQ_\zs{n}} \L_\zs{2}(Q)
\le \l(1+(\varsigma^{*})^2)\,,
\end{equation}
where the parameter $\varsigma^{*}$ is defined \eqref{sec:Ex.01-1}.
\end{proposition}

\proof
To estimate the total weighted noise deviation \eqref{sec:C2EQ}
 note that for $x\in\cB_\zs{p}$
\begin{equation}
\label{ov-I-2Q-1}
\E_\zs{Q}\left(
\sum^{p}_\zs{j=1}\,
x_\zs{j}\,
\wt{\xi}_\zs{j}
\right)^{2}
=\frac{1}{n^{2}}\,\E_\zs{Q}\overline{I}^{2}_\zs{n}(x)\,,
\end{equation}
where
$$
\overline{I}_\zs{n}(x)=\sum^{p}_\zs{j=1}
x_\zs{j}\,\wt{I}_\zs{n}(\psi_\zs{j})
\quad\mbox{and}\quad
\wt{I}_\zs{n}(\psi_\zs{j})=I^{2}_\zs{n}(\psi_\zs{j})-\E_\zs{Q}\,I^{2}_\zs{n}(\psi_\zs{j})
\,.
$$

\noindent
Using
 Proposition \ref{Pr.sec:Stc.5}, the definition of $\sigma_\zs{Q}$ in \eqref{sec:Ex.01-1},
 and  that $\phi^{2}_\zs{*}=2$, we obtain that  there exists a constant $\l>0$ such that for any  $3\le p\le n$
\begin{equation}\label{sec:Pr.16}
\L_\zs{2}(Q)
\le \l \left(1+
\sigma^{2}_\zs{Q}
+
\varpi^{*}_\zs{n} \right)
\end{equation}
where 
$\varpi^{*}_\zs{n}=\sup_\zs{i\ge 3,j\ge 3,\vert i-j\vert\ge 2}\varpi_\zs{n}(\psi_\zs{i},\psi_\zs{j})$
and
\begin{equation*}\label{sec:Ou.10*}
\varpi_\zs{n}(f,g)
=\max_\zs{0\le v+t\le n}\,
\left(
\left|\int^{t}_\zs{0}f(u+v)g(u)\d u\right|
+
\left|\int^{t}_\zs{0} g(u+v) f(u)\d u\right|
\right)\,.
\end{equation*}
Denoting
$$
\upsilon_\zs{j,l}
=
\max_\zs{0\le v+t\le n}\,
\left|
\int^{t}_\zs{0}\psi_\zs{j}(u+v)\psi_\zs{l}(u)\d u
\right|
\,,
$$
we obtain that
$
\varpi^{*}_\zs{n}\le 2\sup_\zs{l\ge 3,j\ge 3,\vert l-j\vert\ge 2}\,\upsilon_\zs{j,l}$.
To estimate the term $\upsilon_\zs{j,l}$
note that for any $3\le j,l\le p$
$$
\upsilon_\zs{j,l}=\max_\zs{0\le v\le 1}\,
\max_\zs{0\le t\le n-v}\,
\left|[t]\V_\zs{j,l}(v)
+
\int^{t}_\zs{[t]}
\psi_\zs{j}(u+v)\,\psi_\zs{l}
(u)\,
\d u
\right|\,
$$
with $\V_\zs{j,l}(v)=\int^{1}_\zs{0}
\psi_\zs{j}(u+v)\,\psi_\zs{l}
(u)\,
\d u$.
Therefore,
\begin{equation}\label{sec:Pr.16-2}
\upsilon_\zs{j,l}
\le
2
+
n\,
\max_\zs{0\le  v\le 1}\,
\left|\V_\zs{j,l}(v)\right|
\,.
\end{equation}
Taking into account the definition of the functions $\psi_\zs{j}$ in \eqref{sec:Imp.3} for the trigonometric basis \eqref{sec:In.5}, we find that
 for $t_\zs{s-1}\le v\le t_\zs{s}$ the function $\V_\zs{j,l}(v)$
can be represented as
$$
\V_\zs{j,l}(v)
=
(t_\zs{s}-v)\,\sum^{p}_\zs{k=1}\, \Trg_\zs{l}(t_\zs{k})
\Trg_\zs{j}(t_\zs{k}+t_\zs{s-1})
+
(v-t_\zs{s-1})\,\sum^{p}_\zs{k=1}\,\Trg_\zs{l}(t_\zs{k})
\Trg_\zs{j}(t_\zs{k}+t_\zs{s})\,.
$$
We recall that for any $u,v\ge 0$ the  functions
\eqref{sec:In.5} satisfy the following equation
$$
\Trg_\zs{j}(u+v)=\k^{*}_\zs{j,1}(v)\Trg_\zs{j-1}(u)
+\k^{*}_\zs{j,2}(v)\Trg_\zs{j}(u)+
\k^{*}_\zs{j,3}(v)\Trg_\zs{j+1}(u)\,,
$$
where $\k^{*}_\zs{j,1}(\cdot)$, $\k^{*}_\zs{j,2}(\cdot)$ and $\k^{*}_\zs{j,3}(\cdot)$  are some trigonometric functions. Using the  property
\eqref{sec:In.5} we obtain that for any $v>0$ for $|l-j|\ge 2$ the function $\V_\zs{j,l}(v)=0$, i.e.  $\varpi^{*}_\zs{n}\le 4$. Using this in
\eqref{sec:Pr.16} and taking into account that $\sigma_\zs{Q}\le \varsigma^{*}$, we obtain the inequality \eqref{Up-Bnd-21_1}.
Hence Proposition \ref{Pr.sec:L2-11}. \endproof

\section{Proofs}
\label{sec:proofs-1}

\subsection{Proof of Proposition \ref{Pr.sec:Ex-CnD-H-2}}

\proof  First, we represent the process \eqref{sec:Ex.1}
as
$$
\xi_\zs{t}=\varrho_1\xi_\zs{t}^{(1)}+\varrho_2\xi_\zs{t}^{(2)}\,,
$$
where $(\xi_\zs{t}^{(1)})_\zs{t\ge 0}$ and $(\xi_\zs{t}^{(2)})_\zs{t\ge 0}$ are independent
Ornstein--Uhlenbeck processes
$$
\d\xi_\zs{t}^{(1)}=a\xi_\zs{t}^{(1)}dt+\d w_\zs{t}
\quad\mbox{and}\quad
\d\xi_\zs{t}^{(2)}=a\xi_\zs{t}^{(2)}dt+\d z_\zs{t}
$$
with $\xi^{(1)}_\zs{0}=\xi^{(2)}_\zs{0}=0$.
Moreover, for any square integrable  function $f$   we denote
\begin{equation}\label{sec:Stc.7}
I_\zs{t}^{(1)}(f)=\int_0^t f(s)d\xi_s^{(1)}
\quad\mbox{and}\quad
I_\zs{t}^{(2)}(f)=\int_0^t f(s)d\xi_s^{(2)}
\,.
\end{equation}
Then the covariance matrix \eqref{sec:Imp.6-1} can be rewritten in the form
$$
\G=\varrho_1^2\,\G_\zs{1}+\varrho_2^2\G_\zs{2},
$$
and the $(i,j)$-th elements of the matrix $\G_\zs{1}$ and $\G_\zs{2}$ are defined as
$\E_\zs{Q}I_n^{(1)}\,(\psi_\zs{i})I_n^{(1)}(\psi_\zs{j})$ and $\E_\zs{Q}I_n^{(2)}\,(\psi_\zs{i})I_n^{(2)}(\psi_\zs{j})$
with $\psi_{j}(t)=\sum_{k=1}^{np}\Trg_\zs{j}(t_k)\Chi_\zs{(t_{k-1},t_k]}(t)$.
Applying the celebrated inequality
of Lidskii and Wieland (see, for example, in
\cite{MarchallOlkin1979}, G.3.a., p.334)
yields
\begin{equation}\label{eq1.24}
\tr\G-\lambda_{\max}(\G)\ge\varrho_1^{2}(
\tr\,\G_\zs{1}-\lambda_\zs{\max}(\G_\zs{1}))
\quad\mbox{a.s.}
\end{equation}
Now,  using
Proposition \ref{Pr.sec:Stc.1}
with $\varrho_\zs{1}=1$ and
$\varrho_\zs{2}=0$,
 we obtain that
\begin{equation}\label{eq2.27}
\tr \G_\zs{1}=
\frac{1}{n}\sum_{j=1}^d\E_\zs{Q}(I_n^{(1)}(\psi_\zs{j}))^2
=d+\sum_{j=1}^d\,b_\zs{j}\,,
\end{equation}
where
$$
b_\zs{j}
=\frac{a}{n}\int_\zs{0}^n\,\psi_\zs{j}(t)\int_\zs{0}^{t}
e^{a(t-s)}\,\psi_\zs{j}(s)(1+e^{2as})\d s \d t\,.
$$
Setting $\Phi(t,v)=\sum_{j=1}^d\psi_\zs{j}(t)\psi_\zs{j}(t-v)$,
we get
\begin{align*}
\tr \G_\zs{1}&=d+\frac{a}{n}\int_0^n
e^{av}\left(\int_\zs{v}^{n}\Phi(t,v)
\left(
1+e^{2a(t-v)}
\right)
\d t\right)
\d v\\[2mm]
&>
d-
2\vert a\vert
\int_0^n
e^{av}
\,\Phi^{*}_\zs{d}(v)
\d v
\quad\mbox{and}\quad
\Phi^{*}_\zs{d}(v)=\max_\zs{t\ge v}\vert \Phi(t,v) \vert\,.
\end{align*}
\noindent
Note, that the function $\Phi^{*}_\zs{d}(\cdot)$ is $1$ - periodic, therefore, in view of the bound \eqref{UbnD-1-TrGn}
$$
\tr \G_\zs{1}>d-2\vert a\vert
\sum^{n}_\zs{k=1}\,e^{a(k-1)}
\,
\int_\zs{0}^1
\,\Phi^{*}_\zs{d}(v)
\d v
>d-
\frac{5+\ln d}{2\check{a}}
\,,
$$
where $\check{a}$ is defined in \eqref{d-0-Df-1}. Therefore,
$\tr\G_\zs{1}>d/2$ for $d\ge d_\zs{0}$.
Now, note that
$$
\lambda_\zs{\max}(\G_\zs{1}) =\sup_\zs{|z|_\zs{d}=1}n^{-1}\E_\zs{Q}\,(I_n^{(1)}(\s))^2
\,,
$$
where $\s(t)=\sum^{d}_\zs{j=1}\,z_\zs{j}\,\psi_\zs{j}(t)$.
Using
Proposition \ref{Pr.sec:Stc.1}, we find that for $\vert z\vert_\zs{d}=1$
$$
\E_\zs{Q}\,(I_n^{(1)}(\s))^2=
2 a\,\int_\zs{0}^{n}\s(t)\,
\check{\varepsilon}_\zs{t}(\s)
\d t+
\int_\zs{0}^{n}\, \s^2(t)\,
\d t
=
2 a\,\int_\zs{0}^{n}\s(t)\,
\check{\varepsilon}_\zs{t}(\s)
\d t+
n
\,.
$$
For $a\le 0$ through the  Cauchy - Bunyakovsky -  Schwarz  inequality
we obtain
\begin{align*}
2
\left\vert
a\,\int_\zs{0}^{n}\s(t)\,
\check{\varepsilon}_\zs{t}(\s)
\d t
\right\vert
&
\le 2
\vert a\vert\,\int_\zs{0}^{n}\,e^{a u}
\check{\varepsilon}_\zs{t}(\s)
\left(
\int^{n}_\zs{u}\,\vert\s(t)\vert \vert\s(t-u)\vert
\d t
\right)
\d u
\vert
\\[3mm]
&
\le 2\int^{n}_\zs{0}\,\s^{2}(t)\,\d t
=2 n
\int^{1}_\zs{0}\,\s^{2}(t)\,\d t
=2 n
\,,
\end{align*}
i.e. $\lambda_\zs{\max}(\G_\zs{1})\le 3$. Hence Proposition \ref{Pr.sec:Ex-CnD-H-2}.
\endproof

\subsection{Proof of Theorem \ref{sec:Mo.Th.1}}

Substituting \eqref{sec:Mo.4} in \eqref{sec:Mo.1} yields for any $\gamma\in\Gamma$
\begin{equation}\label{sec:Mo.11}
\Er\,(\gamma)=J(\gamma)+
2\,\sum^{p}_\zs{j=1}\,\gamma(j)\left(\theta^{*}_\zs{j}\wh{\theta}_\zs{j}-\frac{\wh{\sigma}_\zs{n}}{n}
-\theta^{*}_\zs{j}
\bar{\theta}_\zs{j}\right)
+
\|S\|^2-\rho\wh{P}_\zs{n}(\gamma)\,.
\end{equation}
Now we put $L(\gamma)=\sum^{p}_\zs{j=1}\,\gamma(j)$,
\begin{equation}\label{Main-mart-1-12}
M(\gamma)=\frac{1}{\sqrt{n}}\sum^{p}_\zs{j=1}\,\gamma(j)\bar{\theta}_\zs{j}\xi_\zs{j}\,,
\quad
M_\zs{1}(\gamma)=\frac{1}{\sqrt{n}}\sum^{p}_\zs{j=1}\,\gamma(j)g(j)\wh{\theta}_\zs{j}\xi_\zs{j}\,,
\end{equation}
\begin{equation}
\label{DeFs-B-1-B-2}
B_\zs{1}(\gamma)=\sum^{p}_\zs{j=1}\,\gamma(j)(\E_\zs{Q}\xi_\zs{j}^2-\sigma_\zs{Q})\,
\quad\mbox{and}\quad
B_\zs{2}(\gamma)=\sum^{p}_\zs{j=1}\,\gamma(j)(\xi^{2}_\zs{j}-\E_\zs{Q}\xi_\zs{j}^2)\,.
\end{equation}
Taking into account the definition \eqref{sec:Mo.5}, we can rewrite \eqref{sec:Mo.11} as
\begin{equation*}
\Er\,(\gamma)\,=\,J_\zs{n}(\gamma)+2\frac{\sigma_\zs{Q}-\wh{\sigma}_\zs{n}}{n}L(\gamma)+
2\,M(\gamma)+\frac{2}{n}B_\zs{1}(\gamma)
\end{equation*}
\begin{equation}\label{sec:Mo.12}
+2\sqrt{P_\zs{n}(\gamma)}\frac{B_\zs{2}(\overline{\gamma})}{\sqrt{\sigma_\zs{Q}n}}-2M_\zs{1}(\gamma)\,+\,
\|S\|^2-\rho\wh{P}_\zs{n}(\gamma)
\end{equation}
with $\overline{\gamma}=\gamma/|\gamma |_\zs{n}$. Let $\gamma_0=(\gamma_0(j))_{1\le j \le p}$ be a fixed sequence in $\Gamma$ and $\gamma^*$ be as in \eqref{sec:Mo.7}.
Substituting $\gamma_0$ and $\gamma^*$ in \eqref{sec:Mo.12}, we consider the difference
\begin{align*}
\Er\,(\gamma^*)-\Er\,(\gamma_0)&\leq 2\frac{\sigma_\zs{Q}-\wh{\sigma}_\zs{n}}{n}L(x)+2M(x)
+\frac{2}{n}B_\zs{1}(x)\\[2mm]
&+2\sqrt{P_\zs{n}(\gamma^*)}\frac{B_\zs{2}(\overline{\gamma^*})}{\sqrt{\sigma_\zs{Q}n}}
-2\sqrt{P_\zs{n}(\gamma_0)}\frac{B_\zs{2}(\overline{\gamma_0})}{\sqrt{\sigma_\zs{Q}n}}\\[2mm]
&-2M_\zs{1}(\gamma^*)+2M_\zs{1}(\gamma_0)
-\rho\wh{P}_\zs{n}(\gamma^*)
+\rho\wh{P}_\zs{n}(\gamma_0)\,,
\end{align*}
where $x=\gamma^*-\gamma_0$.
One has that $|L(x)|\leq 2\nu_\zs{*}$ and $|B_\zs{1}(x)|\le \L_\zs{1}(Q)$.
Applying the elementary inequality
\begin{equation}\label{sec:Mo.13}
2|ab|\leq \varepsilon a^2+\varepsilon^{-1} b^2
\end{equation}
with any $\varepsilon>0$, we get
$$
2\sqrt{P_\zs{n}(\gamma)}\frac{B_\zs{2}(\overline{\gamma})}{\sqrt{\sigma_Q n}}
\leq \varepsilon P_\zs{n}(\gamma)+\frac{B_\zs{2}^2(\overline{\gamma})}{\varepsilon\sigma_Q n}
\leq \varepsilon P_\zs{n}(\gamma)+\frac{B^{*}_\zs{2}}{\varepsilon \sigma n}\,,
$$
 where
 $B^{*}_\zs{2}
=
\max_\zs{\gamma\in\Gamma}\,
\left(
B_\zs{2}^2(\overline{\gamma})
+
B_\zs{2}^2(\overline{\gamma}^{2})
\right)$
with $\gamma^2=(\gamma_j^2)_\zs{1\le j\le n}$.
  From  the definition of the function $\L_\zs{2}(Q)$ in  condition $\C_2)$
we obtain
\begin{equation}\label{sec:Mo.13_Ub}
\E_\zs{Q}\,B^{*}_\zs{2}
\le
\sum_\zs{\gamma\in\Gamma}\,
\left(
\E_\zs{Q}B_\zs{2}^2(\overline{\gamma})
+
\E_\zs{Q}B_\zs{2}^2(\overline{\gamma}^{2})
\right)
\le 2\nu \L_\zs{2}(Q)\,.
\end{equation}

Moreover, by the same method we estimate the term $M_\zs{1}$.
Note that
\begin{equation}\label{sec:Mo.13_Ub++c-n}
\sum^{n}_\zs{j=1}\,g^{2}_\zs{\gamma}(j)\,
\wh{\theta}^{2}_\zs{j}
=\c^{2}_\zs{n}
\le \frac{\c^{*}_\zs{n}}{n}
\,,
\end{equation}
where $\c^{*}_\zs{n}=n\max_{\gamma\in \Gamma}\c_\zs{n}^2$.
Therefore,
 through the  Cauchy-Bunyakovsky-Schwarz inequality we can estimate the term $M_\zs{1}(\gamma)$ as
$$
|M_\zs{1}(\gamma)|\le
\frac{|\gamma |_\zs{n}}{\sqrt{n}}\c_\zs{n}
\left(\sum_\zs{j=1}^{n}
\overline{\gamma}^{2}(j)
\,
\xi^{2}_\zs{j}
\right)^{1/2}
=
\frac{|\gamma |_\zs{n}}{\sqrt{n}}\c_\zs{n}
\left(
\sigma_Q+ B_\zs{2}(\overline{\gamma}^{2})
\right)^{1/2}
\,.
$$
So, applying the elementary inequality \eqref{sec:Mo.13} with some arbitrary $\varepsilon>0$, we have
$$
2|M_\zs{1}(\gamma)|
\leq \varepsilon P_\zs{n}(\gamma)+\frac{\c^{*}_\zs{n}}{\varepsilon\sigma_Q n}
(\sigma_Q+B^{*}_\zs{2})\,.
$$

Using the bounds above, one finds
\begin{align*}
\Er\,(\gamma^*)&\leq\Er\,(\gamma_0)+\frac{4\nu_\zs{*} |\wh{\sigma}_\zs{n}-\sigma_Q|}{n}
+2M(x)+\frac{2}{n}\L_\zs{1}(Q)
\\[2mm]
&+\frac{2}{\varepsilon}\,\frac{\c^{*}}{n\sigma_Q}(\sigma_Q+B^{*}_\zs{2})
+\frac{2}{\varepsilon}\,
\frac{B^{*}_\zs{2}}{n\sigma_Q}
\\[2mm]
&+2 \varepsilon P_\zs{n}(\gamma^*)
+2\varepsilon P_\zs{n}(\gamma_0)
-\rho\wh{P}_\zs{n}(\gamma^*)+\rho\wh{P}_\zs{n}(\gamma_0)\,.
\end{align*}
The setting $\varepsilon=\rho/4$
and the estimating where this is possible $\rho$ by $1$
 in this inequality
imply
\begin{align*}
\Er\,(\gamma^*)&\leq\Er\,(\gamma_0)+
 \frac{5\nu_\zs{*}|\wh{\sigma}_\zs{n}-\sigma_Q|}{n}
+2M(x)+\frac{2}{n}\L_\zs{1}(Q)
\\[2mm]
&+\frac{16 (\c^{*}_\zs{n}+1)(\sigma_Q+B^{*}_\zs{2})}{\rho n\sigma_Q}
-\frac{\rho}{2}\wh{P}_\zs{n}(\gamma^*)+\frac{\rho}{2} P_\zs{n}(\gamma_0)
+\rho\wh{P}_\zs{n}(\gamma_0)\,.
\end{align*}
Moreover, taking into account here that
$$
\vert
\wh{P}_\zs{n}(\gamma_0)
-
P_\zs{n}(\gamma_0)
\vert
\le
\frac{\nu_\zs{*}|\wh{\sigma}_\zs{n}-\sigma_Q|}{n}
$$
and $\rho<1/2$,
we obtain that
\begin{multline}\label{sec:Mo.14}
\Er\,(\gamma^*)\leq\Er\,(\gamma_0)+
 \frac{6\nu_\zs{*}|\wh{\sigma}_\zs{n}-\sigma_Q|}{n}
+2M(x)+\frac{2}{n}\L_\zs{1}(Q)
\\[2mm]
+\frac{16 (\c^{*}_\zs{n}+1)(\sigma_Q+B^{*}_\zs{2})}{\rho n\sigma_Q}
-\frac{\rho}{2} P_\zs{n}(\gamma^*)+\frac{3\rho}{2} P_\zs{n}(\gamma_0)\,.
\end{multline}

\noindent
Now we examine the third term in the right-hand side of this inequality. We have
\begin{equation}
\label{upper_bound_M-+01}
2|M(x)|\leq\varepsilon\|S_\zs{x}\|^2+\frac{Z^*}{n\varepsilon}\,,
\end{equation}
where $S_\zs{x}=\sum^{p}_\zs{j=1}\,x_\zs{j}\bar{\theta}_\zs{j}\phi_\zs{j}$
and
$$
Z^*=\sup_{x\in\Gamma_1}\frac{nM^2(x)}{\|S_x\|^2}
\,.
$$
We remind that the set $\Gamma_\zs{1}=\Gamma-\gamma_\zs{0}$.
Using the property
\eqref{sec:In.3},
we get for any fixed $x=(x_\zs{j})_\zs{1\le j\le n}\in\bbr^{n}$
\begin{equation}
\label{M^2+11-00}
\E_\zs{Q}\,M^2(x)=\frac{\E_\zs{Q}\,I^{2}_\zs{n}\left(S_\zs{x}\right)}{n^{2}}
\le
\frac{\varkappa_\zs{Q} \Vert S_\zs{x}\Vert^{2}}{n}
=\frac{\varkappa_\zs{Q}}{n}\,
\sum^{p}_\zs{j=1}\,x^{2}_\zs{j}\,\bar{\theta}^{2}_\zs{j,p}
\end{equation}
and, therefore,
\begin{equation}
\label{up-Z*-00}
\E_\zs{Q}Z^*
\le
\sum_{x\in\Gamma_1}\frac{n M^2(x)}{\|S_x\|^2}
\leq \varkappa_\zs{Q}
\nu
\,.
\end{equation}
Besides,  the norm $\Vert S^{*}_\zs{\gamma^{*}}-S^{*}_\zs{\gamma_\zs{0}}\Vert$ can be estimated from below as
$$
\Vert S^{*}_\zs{\gamma}-S^{*}_\zs{\gamma_\zs{0}}\Vert^{2}
=
\sum^{p}_\zs{j=1}
(x(j)+\beta(j))^{2}\wh{\theta}^{2}_\zs{j}
\ge
\|\wh{S}_\zs{x}\|^2
+2
\sum^{p}_\zs{j=1}
x(j)\beta(j)\wh{\theta}^{2}_\zs{j}\,,
$$
where $\beta(j)=\gamma_0(j)g_j(\gamma_0)-\gamma(j)g_j(\gamma)$.
Then, in view of \eqref{sec:Imp.4}
\begin{align*}
\|S_\zs{x}\|^2&-\|S^{*}_\zs{\gamma}-S^{*}_\zs{\gamma_\zs{0}}\|^2
\le
\|S_\zs{x}\|^2-\|\wh{S}_\zs{x}\|^2
-2\sum^{p}_\zs{j=1}\,x(j)\beta(j)
\wh{\theta}_\zs{j}^2
\\[2mm]
&
\le
-2M(x^{2})-2\sum^{p}_\zs{j=1}\,x(j)\beta(j)
\wh{\theta}_\zs{j}\bar{\theta}_\zs{j}-
\frac{2}{\sqrt{n}}\Upsilon(x)
\,,
\end{align*}
where $\Upsilon(\gamma)=\sum^{p}_\zs{j=1}\,\gamma(j)\beta(j)
\wh{\theta}_\zs{j}\xi_\zs{j}$. The first term in this inequality we can estimate as
$$
2M(x^{2})\le \varepsilon\|S_\zs{x}\|^2+\frac{Z_1^*}{n\varepsilon}
\quad\mbox{and}\quad
Z^*_1=\sup_{x\in\Gamma_1}\frac{n M^2(x^{2})}{\|S_x\|^2}
\,.
$$
Similarly to \eqref{up-Z*-00}, we can estimate the last term as
$$
\E_\zs{Q}Z_1^*\leq \sigma_Q \nu\,.
$$
From this it follows that for any $0<\varepsilon<1$
\begin{equation}\label{upper-bound-000}
\|S_\zs{x}\|^2
\le
\frac{1}{1-\varepsilon}
\left(
\|S^{*}_\zs{\gamma}-S^{*}_\zs{\gamma_\zs{0}}\|^2
+\frac{Z_1^*}{n\varepsilon}
-2\sum^{p}_\zs{j=1}\,x(j)\beta(j)
\wh{\theta}_\zs{j}\bar{\theta}_\zs{j}-
\frac{2\Upsilon(x)}{\sqrt{n}}
\right)
\,.
\end{equation}
The property
 \eqref{sec:Mo.13_Ub++c-n} yields
\begin{equation}
\label{theta-whg-upperb-00}
\sum^{p}_\zs{j=1}\,
\beta^{2}(j)
\wh{\theta}^{2}_\zs{j}
\le
2
\sum^{p}_\zs{j=1}\,g^{2}_\zs{\gamma}(j)\,
\wh{\theta}^{2}_\zs{j}
+
2\sum^{p}_\zs{j=1}\,g^{2}_\zs{\gamma_\zs{0}}(j)\,
\wh{\theta}^{2}_\zs{j}
\le
\frac{4\c^{*}}{\varepsilon n}\,.
\end{equation}
Taking into account that $\vert x(j)\vert\le 1$ and using the inequality \eqref{sec:Mo.13}, we get
 for any $\varepsilon>0$
$$
2\left|\sum^{p}_\zs{j=1}\,x(j)
\beta(j)
\wh{\theta}_\zs{j}\bar{\theta}_\zs{j}\right|\leq\varepsilon\|S_\zs{x}\|^2
+\frac{4\c^{*}}{\varepsilon n}
\,.
$$
To estimate the last term in the right hand of  \eqref{upper-bound-000} we use first
the Cauchy-Bunyakovsky-Schwarz inequality
and then the bound
\eqref{theta-whg-upperb-00}, i.e.
\begin{align*}
\frac{2}{\sqrt{n}}\vert \Upsilon(\gamma)\vert
&\le
\frac{2\nu_\zs{*}}{\sqrt{n}}\left(\sum^{p}_\zs{j=1}\,\beta^{2}(j)
\wh{\theta}^{2}_\zs{j}\right)^{1/2}
\left(
\sum^{p}_\zs{j=1}\,\bar{\gamma}^{2}(j)\,\xi^{2}_\zs{j} \right)^{1/2}\\[2mm]
&
\le
\varepsilon P_\zs{n}(\gamma)
+
\frac{\c^{*}}{n\varepsilon\sigma_Q}
\sum^{p}_\zs{j=1}\,\bar{\gamma}^{2}(j)\,\xi^{2}_\zs{j}
\le
\varepsilon P_\zs{n}(\gamma)
+
\frac{\c^{*}(\sigma_Q+B^{*}_\zs{2})}{n\varepsilon\sigma_Q}
\,.
\end{align*}
Therefore,
$$
\frac{2}{\sqrt{n}}\vert \Upsilon(x)\vert
\le
\frac{2}{\sqrt{n}}\vert \Upsilon(\gamma^{*})\vert
+\frac{2}{\sqrt{n}}\vert \Upsilon(\gamma_\zs{0})\vert
\le
\varepsilon P_\zs{n}(\gamma^{*})
+\varepsilon P_\zs{n}(\gamma_\zs{0})
+
\frac{2\c^{*}(\sigma_Q+B^{*}_\zs{2})}{n\varepsilon\sigma_Q}
\,.
$$
So, using all these bounds in \eqref{upper-bound-000},
we obtain
$$
\|S_\zs{x}\|^{2} \le
\frac{1}{(1-\varepsilon)}\biggl(\frac{Z_1^*}{n\varepsilon}+
\|S_\zs{\gamma^*}^*-S_\zs{\gamma_0}^*\|^2
+\frac{6\c^{*}_\zs{n}(\sigma+B^{*}_\zs{2})}{n\sigma\varepsilon}
+\varepsilon P_\zs{n}(\gamma^*)+\varepsilon P_\zs{n}(\gamma_0)\biggr)\,.
$$
Using in the inequality \eqref{upper_bound_M-+01} this bound and the estimate
$$
\|S_\zs{\gamma^*}^*-S_\zs{\gamma_0}^*\|^2\leq
2(\Er\,(\gamma^*)+\Er\,(\gamma_0))\,,
$$
we have
\begin{align*}
2|M(x)|&\le
\frac{Z^*+Z_1^*}{n(1-\varepsilon)\varepsilon}
+
\frac{2\varepsilon(\Er\,(\gamma^*)+\Er\,(\gamma_0))}{(1-\varepsilon)}
\\[2mm]
&
+
\frac{6\c^{*}_\zs{n}(\sigma_Q+B^{*}_\zs{2})}{n\sigma_Q(1-\varepsilon)}
+\frac{\varepsilon^{2}}{1-\varepsilon}\left( P_\zs{n}(\gamma^*)+ P_\zs{n}(\gamma_0)\right)
\,.
\end{align*}
Choosing here $\varepsilon\le \rho/2<1/2$, we find
\begin{align*}
2|M(x)|&\le
\frac{2(Z^*+Z_1^*)}{n\varepsilon}
+
\frac{2\varepsilon(\Er\,(\gamma^*)+\Er\,(\gamma_0))}{(1-\varepsilon)}
\\[2mm]
&
+
\frac{12\c^{*}_\zs{n}(\sigma_Q+B^{*}_\zs{2})}{n\sigma_Q}
+
\varepsilon\left( P_\zs{n}(\gamma^*)+ P_\zs{n}(\gamma_0)\right)
\,.
\end{align*}
From here and \eqref{sec:Mo.14}, it follows that
\begin{align*}
\Er\,(\gamma^*) &\leq\frac{1+\varepsilon}{1-3\varepsilon}\Er\,(\gamma_0)
+ \frac{6\nu_\zs{*}|\wh{\sigma}_\zs{n}-\sigma_Q|}{n(1-3\varepsilon)}
+\frac{2}{n(1-3\varepsilon)}\L_\zs{1}(Q)
\\[2mm]
&
+\frac{28(1+\c^{*}_\zs{n})(B^{*}_\zs{2}+\sigma_Q)}{\rho(1-3\varepsilon)n\sigma_Q}
+\frac{2(Z^*+Z_1^*)}{n(1-3\varepsilon)}
+\frac{2\rho P_\zs{n}(\gamma_0)}{1-3\varepsilon}.
\end{align*}
Putting here $\varepsilon=\rho/3$
and estimating $(1-\rho)^{-1}$ by $2$ where this is possible,
 we get
\begin{align*}
\Er\,(\gamma^*) &\leq\frac{1+\rho/3}{1-\rho}\Er\,(\gamma_0)
+ \frac{12\nu_\zs{*}|\wh{\sigma}_\zs{n}-\sigma_Q|}{n}
+\frac{4}{n}\L_\zs{1}(Q)
\\[2mm]
&
+\frac{56(1+\c^{*}_\zs{n})(B^{*}_\zs{2}+\sigma_Q)}{\rho n\sigma_Q}
+\frac{4(Z^*+Z_1^*)}{n}
+\frac{2\rho P_\zs{n}(\gamma_0)}{1-\rho}\,.
\end{align*}
Taking the expectation and using the upper bound for $P_\zs{n}(\gamma_0)$  in Lemma \ref{Lem.A.1}
 with $\varepsilon=\rho$ yields
$$
\cR_\zs{Q}(S^*,S)\leq\frac{1+5\rho}{1-\rho}\cR_\zs{Q}(S^*_\zs{\gamma_0},S)
+\frac{\U_\zs{n}}{n\rho}
+
 \frac{12\nu_\zs{*}\E_\zs{Q}|\wh{\sigma}_\zs{n}-\sigma_Q|}{n}
\,,
$$
where
$
\U_\zs{n}=4\L^{*}_\zs{1,n}+
56(1+c^{*}_\zs{n})
(2\L^{*}_\zs{2,n}\nu+1)
+2\c^{*}_\zs{n}$ and the terms $\L^{*}_\zs{1,n}=\sup_\zs{Q\in\cQ_\zs{n}} \L_\zs{1}(Q)$ and  $\L^{*}_\zs{2,n}=\sup_\zs{Q\in\cQ_\zs{n}} \L_\zs{2}(Q)$.
The inequality holds for each $\gamma_0\in\Gamma$, this  implies  Theorem \ref{sec:Mo.Th.1}. \endproof

\subsection{Proof of Proposition \ref{Pr.sigma-est-1}}
\proof
Using the regression scheme \eqref{sec:Imp.4} in  \eqref{sec:Mo.3},
   we obtain that
\begin{equation}
\label{sigma-repr-121}
\wh{\sigma}_\zs{n}=
\frac{n}{p}
\sum_\zs{j=[\sqrt{n}]+1}^n \wh{\t}^2_\zs{j}
=\frac{n}{p}\sum^{p}_\zs{j=[\sqrt{n}]+1}
\bar{\t}^2_\zs{j}
+\frac{2n}{p}
M(\gamma')
+\frac{1}{p}\sum^{p}_\zs{j=[\sqrt{n}]+1}\xi^{2}_\zs{j}
\end{equation}
where the term $M(\gamma)$ is defined in \eqref{Main-mart-1-12},
$\gamma'=(\gamma'(j))_\zs{1\le j\le p}$ and $\gamma'(j)=\Chi_\zs{\{j\ge [\sqrt{n}]+1\}}$.  Note here, that the noise coefficients
$\xi_\zs{j}$ are defined for the trigonometric basis. The first term we estimate trough Lemma A.3 from \cite{KonevPergamenshchikov2015}, i.e.
we have
$$
\frac{n}{p}
\sum^{p}_\zs{j=[\sqrt{n}]+1}
\bar{\t}^2_\zs{j}
\le
\frac{8n}{p}
\left(\int^{1}_\zs{0}\vert \dot{S}(t)\vert \d t\right)^{2}
\sum_\zs{j>[\sqrt{n}]}\,j^{-2}
\le
\frac{8n\Vert \dot{S}\Vert^{2}}{[\sqrt{n}]p}
\,.
$$
Moreover, note now that
$$
M(\gamma')=\frac{1}{n}\,I_\zs{n}(S_\zs{\gamma'})
\quad\mbox{and}\quad
S_\zs{\gamma'}=\sum^{p}_\zs{j=1}\,\gamma'(j)\bar{\t}_\zs{j}\Trg_\zs{j}
\,.
$$
Therefore, the property \eqref{sec:In.3} and condition
$\H_\zs{1})$
 provide that for some constant $\l>0$
$$
\E_\zs{Q,S}\,M^{2}(\gamma')
\le \frac{\varkappa^{*}}{n^{2}}\,\int^{n}_\zs{0}\,S^{2}_\zs{\gamma'}(t)\d t
=
\frac{\varkappa^{*}}{n}
\Vert S_\zs{\gamma'}\Vert^{2}
=\frac{\varkappa^{*}}{n}
\sum^{p}_\zs{j=[\sqrt{n}]+1}
\bar{\t}^2_\zs{j}
\le \l\frac{\Vert\dot{S}\Vert^{2}}{n\sqrt{n}}
\,.
$$
\noindent As to the last term in the representation
\eqref{sigma-repr-121}, we note that
$$
\frac{1}{p}\sum^{p}_\zs{j=[\sqrt{n}]+1}\xi^{2}_\zs{j}=\frac{p-\sqrt{n}}{p}\sigma_\zs{Q}
+
\frac{1}{p}
B_\zs{1}(\gamma')
+
\frac{1}{\sqrt{p}}
B_\zs{2}(\gamma'')\,,
$$
where
$\gamma''=\gamma'/\sqrt{p}$ and the functions $B_\zs{1}(\cdot)$, $B_\zs{2}(\cdot)$ are defined
in \eqref{DeFs-B-1-B-2}.
The using  conditions $\H_\zs{1})$, $\C_\zs{1}$) -- $\C_\zs{3}$)
and $\D$) implies the upper bound \eqref{UpBnd-sigma-101-1}. Hence Proposition \ref{Pr.sigma-est-1}.

\endproof

\section*{Acknowledgements}
This research was supported by RSF, project no 20-61-47043.
The work of S. Pergamenshchikov was partially supported  by
the Ministry of Science and Higher Education the Russian Federation (project no. 1.472.2016/1.4)
and the work of E. Pchelintsev and M. Povzun was partially supported by the Grant of the President of
the Russian Federation (project no. MK-834.2020.9).

\setcounter{section}{0}
\renewcommand{\thesection}{\Alph{section}}

\section{Appendix}\label{sec:Appendix}
\subsection{Stochastic calculus for the non-Gaussian Ornstein--Uhlenbeck--L\'evy processes}\label{sec:OUM}

In this section we study the process \eqref{sec:Ex.1}.
\begin{proposition}\label{Pr.sec:Stc.1}
Let $f$ and $g$ be two nonrandom left continuous $\bbr_\zs{+}\to\bbr$ functions
with the finite right limits.
Then for any $t>0$
\begin{equation}\label{sec:Stc.1}
\E_\zs{Q}\, I_\zs{t}(f)I_\zs{t}(g)=\sigma_\zs{Q}\,
\tau_\zs{t}(f,g)\,,
\end{equation}
where the variance $\sigma_\zs{Q}$ is defined in\eqref{sec:Ex.01-1},
$\tau_\zs{t}(f,g)=
\int^{t}_\zs{0}\left(
f(s)g(s)
+
\check{\varepsilon}_\zs{s}(f)
g(s)+
f(s)
\check{\varepsilon}_\zs{s}(g)
\right)
\,\d s$
and
$$
\check{\varepsilon}_\zs{t}(f)=a\int^{t}_\zs{0}\,e^{a(t-s)}\,f(s)\left(
\frac{1+e^{2as}}{2}
\right)\,\d s
\,.
$$
\end{proposition}
\proof
From the definitions
\eqref{sec:Imp.4} and  \eqref{sec:Ex.1}
we obtain through the Ito formula that
\begin{align}\label{sec:Stc.2}
I_\zs{t}(f)\,I_\zs{t}(g)=
\sigma_\zs{Q}\,
\int^{t}_\zs{0}\,f(s)g(s)\d s+
a\int^{t}_\zs{0}
\Upsilon_\zs{s}(f,g)\,\xi_\zs{t}
\d s
+\M_\zs{t}(f,g)
\,,
\end{align}
where $\Upsilon_\zs{s}(f,g)=f(s) I_\zs{s}(g) + g(s) I_\zs{s}(f)$,
$
\M_\zs{t}(f,g)=\int^{t}_\zs{0}\,\Upsilon_\zs{s-}(f,g)
\,
\d u_\zs{s}
+
\varrho^{2}_\zs{2}\,
\int^{t}_\zs{0}
f(s)\,g(s)\,
\d m_\zs{s}
$
and $m_\zs{t}=x^{2}*(\mu-\wt{\mu})_\zs{t}$.
Moreover, using the Ito formula, we get
\begin{equation}\label{sec:Stc.3}
\E_\zs{Q}\,I^{2}_\zs{t}(1)=
\E_\zs{Q}\, \xi^{2}_\zs{t}
=\sigma_\zs{Q}\,
\frac{e^{2at}-1}{2a}\,.
\end{equation}
Note now that
\begin{align*}
\E_\zs{Q}\,
I^{2}_\zs{t}(f)&=
\E_\zs{Q}
\left(a\int^{t}_\zs{0}\,f(s)\xi_\zs{s}\d s
+
\int^{t}_\zs{0}\,f(s)\,\d u_\zs{s}
\right)^{2}\\[2mm]
&\le
2a^{2}\,
\int^{t}_\zs{0}\,f^{2}(s)\d s
\int^{t}_\zs{0}\,\E\,\xi^{2}_\zs{s}\d s
+2\sigma_\zs{Q}\,
\int^{t}_\zs{0}\,f^{2}(s)\d s\,.
\end{align*}
So, from here
\begin{equation}\label{sec:Stc.4}
\sup_\zs{0\le t\le n}
\E_\zs{Q}\,I^{2}_\zs{t}(f)
\le
2\sigma_\zs{Q}\,(\vert a\vert+
1)\,
\int^{n}_\zs{0}\,f^{2}(s)\d s
<\infty
\,.
\end{equation}
This implies immediately that
$\E_\zs{Q}\,\M_\zs{t}(f,g)=0$. Using this in \eqref{sec:Stc.2} yields
\begin{equation}
\label{sec:Stc.5}
\E_\zs{Q}\,I_\zs{t}(f)\,I_\zs{t}(g)=
\sigma_\zs{Q}\,
\int^{t}_\zs{0}\,f(s)g(s)\d s
+
a\int^{t}_\zs{0}
\left(
f(s)\,\E_\zs{Q} \zeta_\zs{s}(g)
+
g(s)\,
\E_\zs{Q} \zeta_\zs{s}(f)
\right)
\,
\d s
\,,
\end{equation}
where $\zeta_\zs{t}(f)=\xi_\zs{t}\,I_\zs{t}(f)=\,I_\zs{t}(1)\,I_\zs{t}(f)$. Therefore, putting $g=1$ in \eqref{sec:Stc.5}, we have
$$
\E_\zs{Q}\,\zeta_\zs{t}(f)
=
\sigma_\zs{Q}\,
\int^{t}_\zs{0}\,f(s)\d s+
a\int^{t}_\zs{0}
\left(
f(s)
\E_\zs{Q}\,\zeta_\zs{s}(1)
+
\E_\zs{Q}\,\zeta_\zs{s}(f)
\right)
\,
\d s\,.
$$
Taking into account here that $\zeta_\zs{t}(1)=\xi^{2}_\zs{t}$, we find
$$
\E_\zs{Q}\,\zeta_\zs{t}(f)=
\sigma_\zs{Q}\,
\int^{t}\,e^{a(t-s)}f(s)
\frac{1+e^{2as}}{2}
\d s
=\sigma_\zs{Q}\,
\check{\varepsilon}_\zs{t}(f)\,.
$$
Therefore, using this in
 \eqref{sec:Stc.5}, we obtain \eqref{sec:Stc.1}.
\endproof

\begin{corollary}\label{Co.sec:Stc.1}
For any c\`adl\`ag function $f$  from $\L_\zs{2}[0,n]$
\begin{equation}\label{sec:Stc.6}
\E_\zs{Q}\, I^{2}_\zs{n}(f)
\le
2\sigma_\zs{Q}\,
\int^{n}_\zs{0}\,f^{2}(s)
\,\d s\,.
\end{equation}
\end{corollary}
\noindent
This corollary follows directly from  Proposition \ref{Pr.sec:Stc.1} through the Cauchy - Bunyakovsky - Schwarz inequality.

\noindent Now we set
\begin{equation}\label{sec:Stc.6-00}
\wt{I}_\zs{t}(f)=I^{2}_\zs{t}(f)
-\E_\zs{Q}\,I^{2}_\zs{t}(f)
\quad\mbox{and}\quad
V_\zs{t}(f)=\zeta_\zs{t}(f)\,
-
\E_\zs{Q}\,\zeta_\zs{t}(f)
\,.
\end{equation}
Using \eqref{sec:Stc.2} with $f=g$, one has
\begin{equation}\label{sec:Stc.6-01}
\d \wt{I}_\zs{t}(f)
=2af(t)\,V_\zs{t}(f)\d t
+
\d\wt{M}_\zs{t}(f)
\,,
\end{equation}
where
$\wt{M}_\zs{t}(f)=M_\zs{t}(f,f)$.
To study this process we set

\begin{equation}\label{sec:Stc.6-02-0}
\check{\tau}_\zs{t}(f,g)=f(t)g(t)\tau_\zs{t}(1,1)
+
f(t)\tau_\zs{t}(1,g)
+
g(t)\tau_\zs{t}(1,f)
+
\tau_\zs{t}(f,g)
\end{equation}
and
\begin{equation}\label{sec:Stc.6-01-00-1}
A_\zs{t}(f)=
\int^{t}_\zs{0} e^{3a(t-s)}f(s)\upsilon(s)\d s
+2
\,
\sigma_\zs{Q}^2
\,
\int^{t}_\zs{0} e^{3a(t-s)}\check{\varepsilon}_\zs{s}(f)\d s
\,,
\end{equation}
where $\upsilon(s)=a^{2}\E_\zs{Q}\,\wt{\xi}^{2}_\zs{s}+\sigma_\zs{Q}^2
\left( e^{2as}-1 \right) +a\check{\varrho}_\zs{2}$,
 $\wt{\xi}_\zs{s}=\xi^{2}_\zs{s}-\E_\zs{Q}\,\xi^{2}_\zs{s}$ and
 $\check{\varrho}_\zs{2}=\varrho^{4}_\zs{2}\,\Pi(x^{4})$.


\begin{proposition}\label{Pr.sec:Stc.2-1}
For any
 left continuous functions with finite right limits $f$ and $g$
\begin{equation}\label{sec:Stc.6-01-3}
\E_\zs{Q}\, V_\zs{t}(f) V_\zs{t}(g)=
\int^{t}_\zs{0}\,
e^{2a(t-s)}
\,
H_\zs{s}(f,g)
\,\d s
\end{equation}
where
$
H_\zs{t}(f,g)=g(t)A_\zs{t}(f)+
f(t)A_\zs{t}(g)
+\sigma_\zs{Q}^2\check{\tau}_\zs{t}(f,g)
+\check{\varrho}_\zs{2}\,
f(t)g(t)$.
\end{proposition}
\proof Applying again \eqref{sec:Stc.2} with $g=1$
yields
\begin{equation}\label{sec:Stc.6-02}
\d V_\zs{t}(f)
=a
V_\zs{t}(f)\d t
+
a
\,f(t)\,\wt{I}_\zs{t}(1)
\d t
+
\d\,L_\zs{t}(f)
\,,
\end{equation}
where
$
L_\zs{t}(f)
=
\int^{t}_\zs{0}\,
\check{I}_\zs{s-}(f)
\d u_\zs{s}
+
\varrho^{2}_\zs{2}
\int^{t}_\zs{0}\,f(s)
\d m_\zs{s}$ and
$\check{I}_\zs{s}(f)=f(s)\xi_\zs{s}
+I_\zs{s}(f)$.
By the Ito formula we get
\begin{align*}
\d V_\zs{t}(f) V_\zs{t}(g)
&= 2a V_\zs{t}(f) \,V_\zs{t}(g)\d t
+
a\left(
g(t)V_\zs{t}(f)
+
f(t)
V_\zs{t}(g)
\right)\,\wt{I}_\zs{t}(1)\d t
\\[2mm]
&
+\d\,[L(f)\,,\,L(g)]_\zs{t} +
V_\zs{t-}(f) \d L_\zs{t}(g)
+
V_\zs{t-}(g) \d L_\zs{t}(f)\,.
\end{align*}

\noindent
Now, in view of Lemma A.2 from \cite{PchelintsevPergamenshchikov2018}, we find
\begin{align}\nonumber
\d \E_\zs{Q}\,V_\zs{t}(f) V_\zs{t}(g)
&= 2a \E_\zs{Q}\,V_\zs{t}(f) \,V_\zs{t}(g)\d t
+
\left(
g(t)A_\zs{t}(f)
+
f(t)
A_\zs{t}(g)
\right)\,
\d t
\\[2mm]\label{sec:Stc.6-02-0-12-23}
&
+\d\,\E_\zs{Q}\,[L(f)\,,\,L(g)]_\zs{t}\,,
\end{align}
where $A_\zs{t}(f)=a\E_\zs{Q}\,V_\zs{t}(f)\,\wt{I}_\zs{t}(1)=a\E_\zs{Q}\,V_\zs{t}(f)\,V_\zs{t}(1)$.
Note that $\E_\zs{Q}\,\check{I}_\zs{s}(f)\check{I}_\zs{s}(g)=\sigma_\zs{Q}\check{\tau}_\zs{s}(f,g)$
 and
\begin{align*}
\E_\zs{Q}\,[L(f)\,,\,L(g)]_\zs{t}
&=\varrho^{2}_\zs{1}\,\int^{t}_\zs{0}\,
\E_\zs{Q}\,\check{I}_\zs{s}(f)\,
\check{I}_\zs{s}(g)
\,
\d s
+
\E_\zs{Q}\,\sum_\zs{0\le s\le t}\,\Delta\,L_\zs{s}(f)
\Delta\,L_\zs{s}(g)\\[2mm]
&=\sigma_\zs{Q}^2\,\int^{t}_\zs{0}\,\check{\tau}_\zs{s}(f,g)\,\d s
+
\check{\varrho}_\zs{2}
\,\int^{t}_\zs{0}\,f(s)g(s)\d s
\,.
\end{align*}
\noindent
To find  $A_\zs{t}(f)$ we put $g=1$ in \eqref{sec:Stc.6-02-0}.
Note that $A_\zs{t}(1)=\wt{I}^{2}_\zs{t}=\wt{\xi}^{2}_\zs{t}$. Then
$$
\E_\zs{Q}\,V_\zs{t}(f)V_\zs{t}(1)
=\int^{t}_\zs{0}\,e^{3a(t-s)}\,\left(a\,f(s)\,\E_\zs{Q}\,\wt{\xi}^{2}_\zs{s} +
\sigma_\zs{Q}^2\,\check{\tau}_\zs{s}(f,1)
+
\check{\varrho}_\zs{2}\,
f(s)
\right)\d s
\,.
$$
Using here that
\begin{equation}\label{sec:Stc.7-06-1-0}
a\tau_\zs{t}(1,1)=(e^{2at}-1)/2
\quad\mbox{and}\quad
a\,\tau_\zs{t}(1,f)=\check{\varepsilon}_\zs{t}(f)
\,,
\end{equation}
we obtain the representation \eqref{sec:Stc.6-01-3}. Hence Proposition \ref{Pr.sec:Stc.2-1}.
\endproof

\begin{proposition}\label{Pr.sec:Stc.2}
For any
 left continuous function $f$ with finite right limits
\begin{equation}\label{sec:Stc.7-06-1}
\E_\zs{Q}\, \wt{I}_\zs{t}(f)\,\wt{I}_\zs{t}(1)
 =
 \int^{t}_\zs{0}\,
e^{2a(t-s)}
\,
\wt{\varkappa}_\zs{s}(f)
\,
\d s\,,
\end{equation}
where
$
\wt{\varkappa}_\zs{s}(f)
=
2
f(s)\,
A_\zs{s}(f)
\,+
4\sigma_\zs{Q}^2\,
f(s)
\tau_\zs{s}(f,1)
+
\check{\varrho}_\zs{2}\,f^{2}(s)
$.
\end{proposition}
\proof
 Using the Ito formula and Lemma A.2 from \cite{PchelintsevPergamenshchikov2018},
we get for any bounded nonrandom functions $f$ and $g$
\begin{align}\nonumber
\d \E_\zs{Q}\, \wt{I}_\zs{t}(f)\,V_\zs{t}(g)
 &=a \E_\zs{Q}\, \wt{I}_\zs{t}(f)\,V_\zs{t}(g)\,\d t
+2af(t)\E_\zs{Q}\,V_\zs{t}(f)\,V_\zs{t}(g)\,\d t\\[2mm]\label{sec:Stc.7-06-2}
&+
a\,g(t)\,\E_\zs{Q}\,\wt{I}_\zs{t}(f)\,\wt{I}_\zs{t}(1)\d t
+
\d \E_\zs{Q}\,[\wt{M}(f)\,,\,L(g)]_\zs{t}
\,.
\end{align}
Putting here $g=1$ and taking into account that $V_\zs{t}(1)=\wt{I}_\zs{t}(1)$,
we have
\begin{equation*}
\d \E_\zs{Q}\, \wt{I}_\zs{t}(f)\,V_\zs{t}(1)
 =2a \E_\zs{Q}\, \wt{I}_\zs{t}(f)\,V_\zs{t}(1)\,\d t
+2af(t)\E_\zs{Q}\,V_\zs{t}(f)\,V_\zs{t}(1)\,\d t
+
\d \E_\zs{Q}\,[\wt{M}(f)\,,\,L(1)]_\zs{t}
\,.
\end{equation*}
\noindent
Taking into account here that
$\E_\zs{Q}\,[\wt{M}(f)\,,\,L(1)]_\zs{t}=\int^{t}_\zs{0}\,\check{A}_\zs{s}(f)\d s$,  we come to the equality \eqref{sec:Stc.7-06-1}.
\endproof

\noindent
Now, for  integrable $[0,+\infty)\to \bbr$ functions $f$ and $g$ we define the correlation measure as
\begin{equation}\label{sec:Ou.10}
\varpi_\zs{n}(f,g)
=\max_\zs{0\le v+t\le n}\,
\left(
\left|\int^{t}_\zs{0}f(u+v)g(u)\d u\right|
+
\left|\int^{t}_\zs{0} g(u+v) f(u)\d u\right|
\right)
\end{equation}

\begin{proposition}\label{Pr.sec:Stc.3}
Let  $f$ and $g$ be two
 left continuous bounded by $\phi_\zs{*}$ functions with finite right limits, i.e. $\Vert f\Vert_\zs{*,n}\le \phi_\zs{*}$
 and $\Vert g\Vert_\zs{*,n}\le \phi_\zs{*}$, where $\|f\|_\zs{*,n}=\sup_\zs{0\le t\le n}|f(t)|$.
 Then for any $0\le t\le n$
\begin{equation}\label{sec:Stc.7-06-1+1}
\left\vert a\E_\zs{Q}\,\wt{I}_\zs{t}(f)\,V_\zs{t}(g)
\right\vert\le u^{*}_\zs{1}\varpi_\zs{t}(1,g)
+u^{*}_\zs{2}\varpi_\zs{t}(f,g)
+u^{*}_\zs{3}\,,
\end{equation}
where $u^{*}_\zs{1}=4\phi^{2}_\zs{*}(a_\zs{max})\check{\varrho}_\zs{2}+3\sigma_\zs{Q}^2$,
$u^{*}_\zs{2}=44\phi_\zs{*}\sigma_\zs{Q}^2$ and $u^{*}_\zs{3}=3\phi^{3}_\zs{*}\check{\varrho}_\zs{2}$.
\end{proposition}
\proof
First, note that from Ito formula
we find
\begin{multline}\label{sec:Stc.7-06-1+2}
a\E_\zs{Q}\,\wt{I}_\zs{t}(f)\,V_\zs{t}(g)=
a^{2}\int^{t}_\zs{0}\,e^{a(t-s)}\,g(s)\,\left(\E_\zs{Q}\,\wt{I}_\zs{s}(f)\,\wt{I}_\zs{s}(1) \right)\d s
\\[2mm]
+2a^{2}
\int^{t}_\zs{0}\,e^{a(t-s)}\,f(s)\,\left(\E_\zs{Q}\,V_\zs{s}(g)\,V_\zs{s}(f) \right)\d s
+a
\int^{t}_\zs{0}\,e^{a(t-s)}\,
\d \E_\zs{Q}\,[\wt{M}(f)\,,\,L(g)]_\zs{s}
\,.
\end{multline}
Using here
  Lemmas A.4
and  A.6 from \cite{PchelintsevPergamenshchikov2018},
we obtain 
\begin{equation}
\label{sec:Stc.11-06-1}
\vert
a \E_\zs{Q}\,V_\zs{t}(g)\,V_\zs{t}(f)
\vert
\le 15 \sigma_\zs{Q}^2\,\varpi_\zs{t}(f,g)
+\check{\varrho}_\zs{2}\,\Vert f\Vert_\zs{*,t}
\Vert g\Vert_\zs{*,t}
\,.
\end{equation}
One can check directly that
\begin{align*}
\E_\zs{Q}\,[\wt{M}(f)\,,\,L(g)]_\zs{t}&=2\sigma_\zs{Q}\,
\int^{t}_\zs{0}\,g(s)f(s)\,
\left(
\E_\zs{Q}\,I_\zs{s}(f)\,I_\zs{s}(1)
\right)
\d s\\[2mm]
&
+
2\sigma_\zs{Q}\,
\int^{t}_\zs{0}\,f(s)\,
\left(
\E_\zs{Q}\,I_\zs{s}(f)\,I_\zs{s}(g)
\right)
\d s
+\check{\varrho}_\zs{2}
\int^{t}_\zs{0}\,f^{2}(s)\,g(s)
\d s\,.
\end{align*}
From \eqref{sec:Stc.1} we find that
\begin{align*}
\E_\zs{Q}\,[\wt{M}(f)\,,\,L(g)]_\zs{s}&=
2\sigma_\zs{Q}^2\,
\int^{t}_\zs{0}\,
g(s)f(s)\,\tau_\zs{s}(f,1)\,
\d s
\\[2mm]
&
+
2\sigma_\zs{Q}^2\,
\int^{t}_\zs{0}\,
f(s)\,\tau_\zs{s}(f,g)\,
\d s
+\check{\varrho}_\zs{2}
\int^{t}_\zs{0}\,f^{2}(s)\,g(s)
\d s\,.
\end{align*}
Using the last equality in \eqref{sec:Stc.7-06-1-0},
we have
\begin{align*}
a\int^{t}_\zs{0}\,e^{a(t-s)}\,
&\d \E_\zs{Q}\,[\wt{M}(f)\,,\,L(g)]_\zs{s}
=2\sigma_\zs{Q}^2\,\int^{t}_\zs{0}\,e^{a(t-s)}\,
g(s)f(s)\,\check{\varepsilon}_\zs{s}(f)\,
\d s\\[2mm]
&+
2\sigma_\zs{Q}^2\,a\int^{t}_\zs{0}\,e^{a(t-s)}\,
f(s)\,\tau_\zs{s}(f,g)\,
\d s
+
a
\check{\varrho}_\zs{2}
\int^{t}_\zs{0}\,e^{a(t-s)}\,
f^{2}(s)\,g(s)
\d s\,.
\end{align*}
Note now that
$$
\check{\varepsilon}^{\prime}_\zs{t}(f)=a\check{\varepsilon}_\zs{t}(f)
+af(t)(1+e^{2at})/2\,,
$$
i.e.
$\Vert\check{\varepsilon}^{\prime}(f) \Vert_\zs{*,t} \le 2\vert a\vert \Vert f\Vert_\zs{*,t}$. Therefore, in view of  Lemma A.3 from \cite{PchelintsevPergamenshchikov2018}
we get
$$
\left\vert
\int^{t}_\zs{0}\,e^{a(t-s)}\,
g(s)f(s)\,\check{\varepsilon}_\zs{s}(f)\d s
\right\vert
\le 4\varpi_\zs{t}(f,g) \Vert f\Vert_\zs{*,t}\,.
$$
By integrating by parts we can obtain
$\left\vert
\int^{t}_\zs{0}\,
g(s)\,\check{\varepsilon}_\zs{s}(f)\d s
\right\vert\le \varpi_\zs{t}(f,g)$,
and, therefore,
\begin{equation}
\label{sec:Stc.7-00+1}
\vert \tau_\zs{t}(f,g)\vert\le 3\varpi_\zs{t}(f,g)\,.
\end{equation}
So, the last term in \eqref{sec:Stc.7-06-1+2}
can be estimated as
$$
\left\vert
\,a\int^{t}_\zs{0}\,e^{a(t-s)}\,
\d \E_\zs{Q}\,[\wt{M}(f)\,,\,L(g)]_\zs{s}
\right\vert
\le 14\,\sigma_\zs{Q}^2
\varpi_\zs{t}(f,g)\,\Vert f\Vert_\zs{*,t}
+\check{\varrho}_\zs{2}\Vert f\Vert^{2}_\zs{*,t} \Vert g\Vert_\zs{*,t}
\,.
$$
Using  Lemma A.5 from \cite{PchelintsevPergamenshchikov2018} in \eqref{sec:Stc.7-06-1+2}
we come to the bound \eqref{sec:Stc.7-06-1+1}.
\endproof

\begin{proposition}\label{Pr.sec:Stc.4}
Let  $f$ and $g$ be two
 left continuous bounded by $\phi_\zs{*}$ functions with finite right limits, i.e. $\Vert f\Vert_\zs{*,n}\le \phi_\zs{*}$
 and $\Vert g\Vert_\zs{*,n}\le \phi_\zs{*}$. Then for any $t>0$
\begin{equation}\label{sec:Stc.8+1}
\left\vert \E_\zs{Q}\,[\wt{M}(f)\,,\, \wt{M}(g)]_\zs{t}\,
\right\vert\le
\left(
12\sigma_\zs{Q}^2\,\phi^{2}_\zs{*}\varpi_\zs{t}(f,g)
+\phi^{4}_\zs{*}\check{\varrho}_\zs{2}
\right)
t\,.
\end{equation}
\end{proposition}
\proof
From \eqref{Pr.sec:Stc.1}
we obtain that
\begin{equation}\label{sec:Stc.8+1-1}
\E_\zs{Q}\,[\wt{M}(f)\,,\, \wt{M}(g)]_\zs{t}=4\sigma_\zs{Q}^2\,\int^{t}_\zs{0}\,f(s)g(s)\tau_\zs{s}(f,g)\,
\d s
+
\check{\varrho}_\zs{2}\,\int^{t}_\zs{0}\,f^{2}(s)g^{2}(s)\,\d s\,.
\end{equation}
Using here the bound \eqref{sec:Stc.7-00+1}
we get  \eqref{sec:Stc.8+1}. Hence Proposition \ref{Pr.sec:Stc.4}.
\endproof

\begin{corollary}\label{Co.sec:Stc.3-1}
Let  $f$ and $g$ be two
 left continuous bounded by $\phi_\zs{*}$ functions with finite right limits, i.e. $\Vert f\Vert_\zs{*,n}\le \phi_\zs{*}$
 and $\Vert g\Vert_\zs{*,n}\le \phi_\zs{*}$. Then for any $t>0$
\begin{equation}\label{sec:Stc.15-06-1+1}
\left\vert\E_\zs{Q}\,\wt{I}_\zs{t}(f)\wt{I}_\zs{t}(g)\right\vert\,\le\,\left(v^{*}_\zs{1}(\varpi_\zs{t}(1,f)+\varpi_\zs{t}(1,g))+v^{*}_\zs{2}\varpi_\zs{t}(f,g)+v^{*}_\zs{3}\right)\,t\,,
\end{equation}
where $v^{*}_\zs{1}=8\phi^{3}_\zs{*}a_\zs{max}\check{\varrho}_\zs{2}+6\sigma_\zs{Q}^2$,
$v^{*}_\zs{2}=100\phi^2_\zs{*}\sigma_\zs{Q}^2$ and $v^{*}_\zs{3}=13\phi^{4}_\zs{*}\check{\varrho}_\zs{2}$.
\end{corollary}
\proof
From \eqref{sec:Stc.6-01} by the Ito formula one finds for $t\ge 0$
$$
\E_\zs{Q}\wt{I}_\zs{t}(f)\wt{I}_\zs{t}(g)=\E_\zs{Q}[\wt{M}(f),\wt{M}(g)]_\zs{t}
+
2a\int^{t}_\zs{0}\left(f(s)\E_\zs{Q}\wt{I}_\zs{s}(g)V_\zs{s}(f)+g(s)\E_\zs{Q}\wt{I}_\zs{s}(f) V_\zs{s}(g)\right)\d s\,.
$$
\noindent
 Using here Proposition \ref{Pr.sec:Stc.3} and Proposition \ref{Pr.sec:Stc.4} we come to desire result.

\endproof

\noindent
Now for square integrable functions $(\phi_\zs{j})_\zs{1\le j\le n}$
and $x=(x_\zs{1},\ldots, x_\zs{p})\in [0,1]^{p}$
 we set
\begin{equation}\label{sec:Stc.9-00+1}
\overline{I}_\zs{p,n}(x)=
\sum^{p}_\zs{j=1}\,x_\zs{j}\,\wt{I}_\zs{n}(\phi_\zs{j})\,.
\end{equation}

\noindent
We need to study this total deviation.

\begin{proposition}\label{Pr.sec:Stc.5}
Assume that $\phi_\zs{1}\equiv 1$ and the functions $(\phi_\zs{j})_\zs{2\le j\le n}$ are bounded, i.e.
$\max_\zs{2\le j\le p}\Vert \phi_\zs{j}\Vert_\zs{*,n}\le \phi_\zs{*}$.
 Then, for any $p\ge 3$ and $n\ge 1$
\begin{equation}\label{sec:Stc.10-23-1}
\sup_\zs{x\in\cB_\zs{p}}\,
\E_\zs{Q}\,\overline{I}^{2}_\zs{p,n}(x)
\le 4  n^{2}\,(4\phi_\zs{*}^2v^{*}_\zs{1}+5\phi_\zs{*}^2v^{*}_\zs{2}+ 2v^{*}_\zs{3})
+
2np\left(\varpi^{*}_\zs{n}(2v^{*}_\zs{1}+v^{*}_\zs{2})+v^{*}_\zs{3} \right),
\end{equation}
where $\varpi^{*}_\zs{n}=\sup_\zs{i\ge 3,j\ge 3,\vert i-j\vert\ge 2}\varpi_\zs{n}(\phi_\zs{i},\phi_\zs{j})$
and
the ball $\cB_\zs{p}$ is defined in
\eqref{sec:C2EQ}.
\end{proposition}

\proof
We
represent the sum \eqref{sec:Stc.9-00+1} as
$$
\overline{I}_\zs{p,n}(x)=
\sum^{p}_\zs{j=1}\,x_\zs{j}\,\wt{I}_\zs{n}(\phi_\zs{j})=
x_\zs{1}\wt{I}_\zs{n}(\phi_\zs{1})+x_\zs{2}
\wt{I}_\zs{n}(\phi_\zs{2})
+
J_\zs{p,n}\,,
$$
where  $J_\zs{p,n}=\sum^n_\zs{j=
3}\,x_\zs{j}\,\wt{I}_\zs{n}(\phi_\zs{j})$. In view of $x^{2}_\zs{1}+x^{2}_\zs{2}\le 1$, one has
\begin{equation*}\label{sec:Pr.3-J12}
\E_\zs{Q}
\overline{I}^{2}_\zs{n}(x)
\,
\le
2\E_\zs{Q}\,\wt{I}^{2}_\zs{n}(\phi_\zs{1}) +
2\E_\zs{Q}\,\wt{I}^{2}_\zs{n}(\phi_\zs{2})
+
2\E_\zs{Q}\,J^{2}_\zs{p,n}\,.
\end{equation*}
\noindent
Using
Corollary~\ref{Co.sec:Stc.3-1}
and taking into account that
$\varpi_\zs{n}(\phi_\zs{i},\phi_\zs{j})\le 2\phi_\zs{*}^2n$ if $|i-j|\le 1$, we get
\begin{equation}\label{sec:Pr.3-J12}
\E_\zs{Q}
\overline{I}^{2}_\zs{n}(x)
\,
\le\,
 4n^2 \left(
4\phi_\zs{*}^2v^{*}_\zs{1}+2\phi_\zs{*}^2v^{*}_\zs{2}+v^{*}_\zs{3}
\right)
+
2\E_\zs{Q}\,J^{2}_\zs{p,n}
\end{equation}
\noindent
and
\begin{equation}\label{sec:Pr.4-J2}
\E_\zs{Q}\,J^{2}_\zs{p,n}=\sum^{p}_\zs{i,j=
3}\,x_\zs{i}x_\zs{j}\,\E_\zs{Q}\wt{I}_\zs{n}(\phi_\zs{i})\wt{I}_\zs{n}(\phi_\zs{j})
\le n
\sum^{p}_\zs{i,j= 3} |x_\zs{i}| |x_\zs{j}|
\wt{\kappa}_\zs{i,j}\,,
\end{equation}
where
$\wt{\kappa}_\zs{i,j}=v^{*}_\zs{1}(\varpi_\zs{n}(1,\phi_\zs{i})+\varpi_\zs{n}(1,\phi_\zs{j}))
+v^{*}_\zs{2}\varpi_\zs{n}(\phi_\zs{i},\phi_\zs{j})+v^{*}_\zs{3}$.
Note here that
for any $i,j\ge 3$ we can estimate the correlation coefficient as
$\varpi_\zs{n}(\phi_\zs{i},\phi_\zs{j})\le 2\phi_\zs{*}^2n\Chi_\zs{\{|i-j|\le 1\}}+
\varpi^{*}_\zs{n}\Chi_\zs{\{|i-j|\ge 2\}}$. By making use of this estimate in
\eqref{sec:Pr.4-J2} and taking into account that
$$
\sum^{p}_\zs{i,j= 3} \Chi_\zs{\{|i-j|\le 1\}}|x_\zs{i}| |x_\zs{j}|\le 3
\quad\mbox{and}\quad
\sum^{p}_\zs{i,j= 3} |x_\zs{i}| |x_\zs{j}|\le p
\,,
$$
 one gets
$$
\sum^{p}_\zs{i,j= 3} |x_\zs{i}| |x_\zs{j}|
\wt{\kappa}_\zs{i,j}
\,
\le
p\left(\varpi^*_\zs{n}(2v^{*}_\zs{1}+v^{*}_\zs{2})+v^{*}_\zs{3}\right)
+ 6\phi_\zs{*}^2nv^{*}_\zs{2}.
$$
From here and the bounds \eqref{sec:Pr.3-J12}--\eqref{sec:Pr.4-J2} we obtain \eqref{sec:Stc.10-23-1}.
Hence Proposition \ref{Pr.sec:Stc.5}. \endproof

\subsection{Property of the trigonometric basis}

\begin{lemma}\label{Le.sec:A.0}
The
 trigonometric basis \eqref{sec:In.5} satisfies the following inequality
\begin{equation}
\label{UbnD-1-TrGn}
\sup_\zs{d\ge 1}\,
\left(
\int^{1}_\zs{0}\,
\max_\zs{t\ge v}
\vert\Phi_\zs{d}(t,v)\vert
-\ln d
\right)
\,
\le 5\,,
\end{equation}
where $\Phi_\zs{d}(t,v)=\sum^{d}_\zs{j=1}\,\Trg_\zs{j}(t)\, \Trg_\zs{j}(t-v)$.
\end{lemma}

\proof
Note that for any $d\ge 3$ and $N=[d/2]$
 the sum $\Phi_\zs{d}(t,v)$ can be represented as
$$
\Phi_\zs{d}(t,v)
=1+2\sum_{j=1}^N\cos(2\pi jv)
-2\sin(2\pi Nt)\sin(2\pi N(t-v))\,\Chi_\zs{\{d=2N\}}\,.
$$
Therefore,
$$
\max_\zs{t\ge v}\,\vert\Phi_\zs{d}(t,v)\vert
\le 2+\vert 1+2\sum_{j=1}^N\cos(2\pi jv)\vert
=2+
\left\vert
\frac{\sin(\pi (2N+1)v)}{\sin(\pi v)}
\right\vert
\,.
$$
Using that $\vert 1+2\sum_{j=1}^N\cos(2\pi jv)\vert\le 2N+1\le d$,
we obtain
 for  $0<\delta<1/2$
$$
\int^{1}_\zs{0}\,
\max_\zs{t\ge v}\,\vert\Phi_\zs{d}(t,v)\vert \d v
\le
2+2\delta (d+1)+2
\int^{1/2}_\zs{\delta}\,\frac{1}{\sin(\pi v)}\d v
\,.
$$
In view of $\sin(\pi v)\ge 2 v$ for any $0<v<1/2$ we get
$$
\int^{1}_\zs{0}\,
\max_\zs{t\ge v}\,\vert\Phi_\zs{d}(t,v)\vert
\,\d v\le 4+2\delta d
-\ln(2\delta)\,.
$$
The minimizing  this upper bound on $0<\delta<1$ implies \eqref{UbnD-1-TrGn}.
Hence Lemma \ref{Le.sec:A.0}. \endproof

\subsection{Property of Penalty term}

\begin{lemma}\label{Lem.A.1}
Assume that conditions $\H_\zs{2})$, $\H_\zs{3})$ and $\C_\zs{1})$ hold. Then,
for any $n\geq 1$, $\gamma\in\Gamma$ and $0<\varepsilon<1$ the penalty term \eqref{sec:Mo.9}
is bounded from above as
\begin{equation}
\label{penalty-00}
P_\zs{n}(\gamma)\leq\frac{\E_\zs{Q,S}\,\Er(\gamma)}{1-\varepsilon}+\frac{\c^{*}_\zs{n}+\L^{*}_\zs{1,n}}{n\varepsilon(1-\varepsilon)}
+
\frac{\L^{*}_\zs{1,n}}{n}
\,,
\end{equation}
where the terms $\c^{*}_\zs{n}$ and $\L^{*}_\zs{1,n}$ are defined in \eqref{Unif-c-1}  and \eqref{sec:C1EQ} respectively.
\end{lemma}

\proof
By the definition of $\Er(\gamma)$ one has
\begin{align*}
\Er(\gamma)&=\sum^{p}_\zs{j=1}\,(\gamma(j)\theta^{*}_\zs{j}-\theta_j)^2
=\sum^{p}_\zs{j=1}\,\left(\gamma(j)(\theta^{*}_\zs{j}-\theta_j)+(\gamma(j)-1)\theta_j\right)^2 \\[2mm]
&
\ge
\sum^{p}_\zs{j=1}\,\gamma^{2}(j)(\theta^{*}_\zs{j}-\theta_j)^2+
2\sum^{p}_\zs{j=1}\,\gamma(j)(\gamma(j)-1)\theta_j(\theta^{*}_\zs{j}-\theta_j).
\end{align*}
Taking into account the condition $\H_\zs{3})$ and the definition
\eqref{sec:Imp.12}, we find
$$
\sum^{p}_\zs{j=1}\,\gamma(j)(\gamma(j)-1)\theta_j(\theta^{*}_\zs{j}-\theta_j)
=
\sum^{p}_\zs{j=1}\,\gamma(j)(\gamma(j)-1)\theta_j(\wh{\theta}_\zs{j}-\theta_j)
\,,
$$
i.e.
$\E_\zs{Q,S}\,\sum^{p}_\zs{j=1}\,\gamma(j)(\gamma(j)-1)\theta_j(\theta^{*}_\zs{j}-\theta_j)=0$. So, in view of  the inequality  \eqref{sec:Mo.13} we get
\begin{align*}
\E_\zs{Q,S}\,\Er(\gamma)&\geq
\sum^{p}_\zs{j=1}\,\gamma^{2}(j)\E_\zs{Q,S}\,(\theta^{*}_\zs{j}-\theta_j)^2=
\sum^{p}_\zs{j=1}\,\gamma^{2}(j)\E_\zs{Q,S}\,\left(\frac{\xi_\zs{j}}{\sqrt{n}}-g(j)\wh{\theta}_\zs{j}\right)^2\\[2mm]&
\ge \frac{1}{n}\,\sum^{p}_\zs{j=1}\,\gamma^{2}(j)
\E_\zs{Q,S}\,\xi^{2}_\zs{j}
-\frac{2}{\sqrt{n}}\E_\zs{Q,S}\,\sum^{p}_\zs{j=1}\,\gamma^{2}(j)g(j)\wh{\theta}_\zs{j}
\xi_\zs{j}
\\[2mm]
&\ge \frac{1-\varepsilon}{n}\,\sum^{p}_\zs{j=1}\,\gamma^{2}(j)
\E_\zs{Q,S}\,\xi^{2}_\zs{j}
-\frac{1}{\varepsilon}\E_\zs{Q,S}\,\sum^{p}_\zs{j=1}\,\gamma^{2}(j) g^{2}(j)\wh{\theta}^{2}_\zs{j}\,.
\end{align*}
The using here here the upper bound \eqref{sec:C1EQ} and the definition of $g(j)$ in
\eqref{sec:Imp.12} yields
\begin{align*}
\E_\zs{Q,S}\,\Er(\gamma)&\ge
(1-\varepsilon)P_\zs{n}(\gamma)-
\frac{(1-\varepsilon)\L^{*}_\zs{1,n}}{n}
-\frac{1}{\varepsilon}\E_\zs{Q,S}\,\sum^{p}_\zs{j=1}\,g^{2}(j)\wh{\theta}^{2}_\zs{j}
\,.
\end{align*}
Now,
the inequality \eqref{sec:Mo.13_Ub++c-n}
 implies the bound \eqref{penalty-00}. Hence  Lemma \ref{Lem.A.1}.
\endproof

\bigskip

\end{document}